\documentclass[12pt]{amsart}
\usepackage[nobysame,abbrev,alphabetic]{amsrefs}
\usepackage{amssymb}
\usepackage[all]{xy}
\usepackage{stmaryrd}

\DeclareMathOperator{\Char}{char}
\DeclareMathOperator{\cd}{cd}

\DeclareMathOperator{\Coker}{Coker}

\DeclareMathOperator{\Gal}{Gal}
\DeclareMathOperator{\Hom}{Hom}
\DeclareMathOperator{\id}{id}
\DeclareMathOperator{\Img}{Im}

\DeclareMathOperator{\Ker}{Ker}

\DeclareMathOperator{\res}{res}

\DeclareMathOperator{\trg}{trg}

\DeclareFontFamily{U}{wncy}{}
\DeclareFontShape{U}{wncy}{m}{n}{<->wncyr10}{}
\DeclareSymbolFont{mcy}{U}{wncy}{m}{n}
\DeclareMathSymbol{\Sha}{\mathord}{mcy}{"58}
\DeclareMathSymbol{\sha}{\mathord}{mcy}{"78}

\begin{document}

\newtheorem{thm}{Theorem}[section]
\newtheorem{cor}[thm]{Corollary}
\newtheorem{lem}[thm]{Lemma}
\newtheorem{prop}[thm]{Proposition}
\newtheorem{defin}[thm]{Definition}
\newtheorem{exam}[thm]{Example}
\newtheorem{examples}[thm]{Examples}
\newtheorem{rem}[thm]{Remark}
\newtheorem{case}{\sl Case}
\newtheorem{claim}{Claim}
\newtheorem{prt}{Part}
\newtheorem*{mainthm}{Main Theorem}
\newtheorem*{thmA}{Theorem A}
\newtheorem*{transfer theorem}{The Transfer Theorem}

\newtheorem{question}[thm]{Question}
\newtheorem*{notation}{Notation}
\swapnumbers
\newtheorem{rems}[thm]{Remarks}
\newtheorem*{acknowledgment}{Acknowledgment}

\newtheorem{questions}[thm]{Questions}
\numberwithin{equation}{section}

\newcommand{\ab}{\mathrm{ab}}
\newcommand{\Ann}{\mathrm{Ann}}
\newcommand{\Bock}{\mathrm{Bock}}
\newcommand{\dec}{\mathrm{dec}}
\newcommand{\diam}{\mathrm{diam}}
\newcommand{\dirlim}{\varinjlim}
\newcommand{\discup}{\ \ensuremath{\mathaccent\cdot\cup}}
\newcommand{\divis}{\mathrm{div}}
\newcommand{\et}{\mathrm{et}}
\newcommand{\gr}{\mathrm{gr}}
\newcommand{\nek}{,\ldots,}
\newcommand{\ind}{\hbox{ind}}
\newcommand{\Ind}{\mathrm{Ind}}
\newcommand{\inv}{^{-1}}
\newcommand{\isom}{\cong}
\newcommand{\Massey}{\mathrm{Massey}}
\newcommand{\ndiv}{\hbox{$\,\not|\,$}}
\newcommand{\nil}{\mathrm{nil}}
\newcommand{\pr}{\mathrm{pr}}
\newcommand{\sep}{\mathrm{sep}}
\newcommand{\sh}{\mathrm{sh}}
\newcommand{\tagg}{^{''}}
\newcommand{\tensor}{\otimes}
\newcommand{\alp}{\alpha}
\newcommand{\gam}{\gamma}
\newcommand{\Gam}{\Gamma}
\newcommand{\del}{\delta}
\newcommand{\Del}{\Delta}
\newcommand{\eps}{\epsilon}
\newcommand{\lam}{\lambda}
\newcommand{\Lam}{\Lambda}
\newcommand{\sig}{\sigma}
\newcommand{\Sig}{\Sigma}
\newcommand{\bfA}{\mathbf{A}}
\newcommand{\bfB}{\mathbf{B}}
\newcommand{\bfC}{\mathbf{C}}
\newcommand{\bfF}{\mathbf{F}}
\newcommand{\bfP}{\mathbf{P}}
\newcommand{\bfQ}{\mathbf{Q}}
\newcommand{\bfR}{\mathbf{R}}
\newcommand{\bfS}{\mathbf{S}}
\newcommand{\bfT}{\mathbf{T}}
\newcommand{\bfZ}{\mathbf{Z}}
\newcommand{\dbA}{\mathbb{A}}
\newcommand{\dbC}{\mathbb{C}}
\newcommand{\dbF}{\mathbb{F}}
\newcommand{\dbG}{\mathbb{G}}
\newcommand{\dbK}{\mathbb{K}}
\newcommand{\dbN}{\mathbb{N}}
\newcommand{\dbP}{\mathbb{P}}
\newcommand{\dbQ}{\mathbb{Q}}
\newcommand{\dbR}{\mathbb{R}}
\newcommand{\dbU}{\mathbb{U}}
\newcommand{\dbV}{\mathbb{V}}
\newcommand{\dbZ}{\mathbb{Z}}
\newcommand{\grf}{\mathfrak{f}}
\newcommand{\gra}{\mathfrak{a}}
\newcommand{\grA}{\mathfrak{A}}
\newcommand{\grB}{\mathfrak{B}}
\newcommand{\grh}{\mathfrak{h}}
\newcommand{\grH}{\mathfrak{H}}
\newcommand{\grI}{\mathfrak{I}}
\newcommand{\grL}{\mathfrak{L}}
\newcommand{\grm}{\mathfrak{m}}
\newcommand{\grp}{\mathfrak{p}}
\newcommand{\grq}{\mathfrak{q}}
\newcommand{\grr}{\mathfrak{r}}
\newcommand{\grR}{\mathfrak{R}}
\newcommand{\grU}{\mathfrak{U}}
\newcommand{\grZ}{\mathfrak{Z}}
\newcommand{\calA}{\mathcal{A}}
\newcommand{\calB}{\mathcal{B}}
\newcommand{\calC}{\mathcal{C}}
\newcommand{\calE}{\mathcal{E}}
\newcommand{\calG}{\mathcal{G}}
\newcommand{\calH}{\mathcal{H}}
\newcommand{\calJ}{\mathcal{J}}
\newcommand{\calK}{\mathcal{K}}
\newcommand{\calL}{\mathcal{L}}
\newcommand{\calR}{\mathcal{R}}
\newcommand{\calW}{\mathcal{W}}
\newcommand{\calU}{\mathcal{U}}
\newcommand{\calV}{\mathcal{V}}
\newcommand{\calZ}{\mathcal{Z}}

\title[The kernel generating condition]{The kernel generating condition and absolute Galois groups}
\author{Ido Efrat}
\address{Earl Katz Family Chair in Pure Mathematics,
Department of Mathematics, Ben-Gurion University of the Negev, P.O.\  Box 653, Be'er-Sheva 8410501, Israel}
\email{efrat@bgu.ac.il}



\maketitle

\begin{abstract}
For a list $\calL$ of finite groups and for a profinite group $G$, we consider the intersection $T(G)$ of all open normal subgroups $N$ of $G$ with $G/N$ in $\calL$.
We give a cohomological characterization of the epimorphisms $\pi\colon S\to G$ of profinite groups (satisfying some additional requirements) such that $\pi[T(S)]=T(G)$.
For $p$ prime, this is used to describe cohomologically the profinite groups $G$ whose $n$th term $G_{(n,p)}$ (resp., $G^{(n,p)}$) in the $p$-Zassenhaus filtration (resp., lower $p$-central filtration) is an intersection of this form.
When $G=G_F$ is the absolute Galois group of a field $F$ containing a root of unity of  order $p$, we recover as special cases results  by Min\'a\v c, Spira and the author, describing $G_{(3,p)}$ and $G^{(3,p)}$ as $T(G)$ for appropriate lists $\calL$.
\end{abstract}

\section{Introduction}
Given a field $F$, let $G_F=\Gal(F^{\rm sep}/F)$ be its absolute Galois group, considered as a profinite group.
It is a major open problem in modern Galois theory to know which profinite groups $G$ are realizable as $G_F$ for some field $F$.
Throughout this paper we fix a prime number $p$ and assume that the field $F$ contains a root of unity of order $p$, and in particular $\Char\,F\neq p$.
Among the few known group-theoretic restrictions on the structure of $G=G_F$ are several ``Intersection Properties" related to standard filtrations of $G$ by characteristic subgroups:
The \textsl{lower $p$-central filtration $G^{(n,p)}$}, and the \textsl{$p$-Zassenhaus filtration} $G_{(n,p)}$, $n=1,2,\ldots \ $.
We recall that they are defined inductively (in the profinite sense) by
\[
\begin{split}
G^{(1,p)}=G&, \quad G^{(n+1,p)}=(G^{(n,p)})^{p}[G,G^{(n,p)}], \\
G_{(1,p)}=G&, \quad  G_{(n,p)}=(G_{(\lceil n/p\rceil,p)})^p\prod_{i+j=n}[G_{(i,p)},G_{(j,p)}]  \hbox{ for } n\geq2.
\end{split}
\]

By results by Min\'a\v c, Spira and the author, for an absolute Galois group $G=G_F$ as above, some of these filtration subgroups can be expressed as the intersection of the open subgroups $N$ of $G$ such that $G/N$ belongs to a certain list $\calL$ of finite groups.
Namely:
\begin{enumerate}
\item[(i)]
When $p=2$,  $G^{(3,2)}=\bigcap\bigl\{N\trianglelefteq  G\ \bigm|\ G/N\isom \dbZ/2,\ \dbZ/4,\ D_4\bigr\}$ \cite{MinacSpira96}*{Cor.\ 2.18};
\item[(ii)]
When $p>2$, $G^{(3,p)}=\bigcap\bigl\{N\trianglelefteq  G\ \bigm|\ G/N\isom \dbZ/p^2,\ M_{p^3}\bigr\}$ \cite{EfratMinac11};
\item[(iii)]
For every $p$, $G_{(3,p)}=\bigcap\bigl\{N\trianglelefteq  G\ \bigm|\ G/N\isom \dbZ/p,\ \dbU_2(\dbZ/p)\bigr\}$ \cite{EfratMinac17}*{Th.\ D}.
\end{enumerate}

Here $D_4$ is the dihedral group of order $8$,  $M_{p^3}$ is the extra-special group of order $p^3$ and exponent $p^2$  (see \S\ref{section on the lower central filtration for odd prime}), and $\dbU_n(R)$ denotes the group of all $(n+1)\times(n+1)$ unipotent upper-triangular matrices over a ring $R$.

Furthermore, when $G$ is a \textsl{free} profinite group and for arbitrary $n\geq2$, one has
\begin{enumerate}
\item[(iv)]
$G_{(n,p)}=\bigcap\bigl\{N\leq  G\ |\ G/N\leq \dbU_n(\dbZ/p)\bigr\}$ \cite{Efrat14a}*{Th.\ A}.
\end{enumerate}
Min\'a\v c and T\^an  \cite{MinacTan15} conjectured that (iv) should hold more generally for every absolute Galois group $G=G_F$ as above.

It should be stressed that the intersection theorems (i)--(iii) do not hold for general profinite groups.
Their proofs rely on deep cohomological properties of absolute Galois groups, in particular the behavior of special elements in $H^2(G_F,\dbZ/p)$, such as cup products, Bockstein elements, and elements of $n$-fold Massey products.

While the known proofs of these intersection theorems are rather different one from each other, they share a common general structure:
First one proves the equality
\begin{equation}
\label{desired intersection property}
G^{(n,p)}=\bigcap\bigl\{N\trianglelefteq G\ \bigm|\ G/N\in\calL\bigr\} \ \ \bigl(\hbox{or } G_{(n,p)}=\bigcap\bigl\{N\trianglelefteq G\ \bigm|\ G/N\in\calL\bigr\}   \bigr)
\end{equation}
when $G$ is a free profinite group.
Then, for a more general profinite group $G$ (such as $G_F$), one takes a profinite presentation, i.e., a continuous epimorphism $\pi\colon S\to G$, where $S$ is a free profinite group, and transfers the equality (\ref{desired intersection property}) from $S$ to $G$.

The first part is purely group-theoretic, and is usually proved using \textsl{Magnus theory}, i.e., by viewing the elements of  $G=S$ as formal power series.
A general machinery to obtain such results in the free profinite case is given in \cite{Efrat14b}  (see also \cite{ChapmanEfrat16}).

In the current paper we focus on the second part, and prove a general cohomological \textsl{Transfer Principle}, which allows one to transfer equalities as in (\ref{desired intersection property}) from a group $S$ to its quotient $G$.
To explain it, we note that for an epimorphism $\pi\colon S\to G$, we clearly have $\pi[S^{(n,p)}]=G^{(n,p)}$ (resp., $\pi[S_{(n,p)}]=G_{(n,p)}$).
However, in general, one does not have that
\begin{equation}
\label{intersection commutes with pi}
\pi\biggl[\bigcap\bigl\{N\trianglelefteq S\ \bigm|\ S/N\in\calL\bigr\}\biggr]=\bigcap\bigl\{M\trianglelefteq G\ \bigm|\ G/M\in\calL\bigr\}.
\end{equation}
Our transfer principle gives a cohomological condition on special elements of $H^2(G,\dbZ/p)$ which is equivalent to the validity of (\ref{intersection commutes with pi}).
From this principle, one can deduce as special cases the intersection theorems (i)--(iii) of \cite{MinacSpira96},  \cite{EfratMinac11} and  \cite{EfratMinac17}, respectively,  as well as several main results of \cite{Efrat14a}.

For the transfer principle, we consider intersections $T(G)$ and $\bar T(G)$ of certain closed normal subgroups of $G$ (see below) with $T(G)\leq\bar T(G)$, and an associated collection $\Pi$ of ``special cohomology elements" in $H^2(G/\bar T(G))$, where $H^i(\bar G)=H^i(\bar G,\dbZ/p)$ denotes the $i$th profinite cohomology group of the profinite group $\bar G$ with the trivial action on $\dbZ/p$.
We further consider the inflation map $\inf\colon H^2(G/\bar T(G))\to H^2(G)$.
For a subset $\Pi_0$ of  $H^2(G/\bar T(G))$ let $\langle\Pi_0\rangle$ be the subgroup it generates.

\begin{transfer theorem}
Let $\pi\colon S\to G$ be a continuous homomorphism of profinite groups, where $\Ker(\pi)\leq\bar T(S)$.
The following conditions are equivalent:
\begin{enumerate}
\item[(a)]
$\pi[T(S)]=T[G]$;
\item[(b)]
$\langle\Pi\cap\Ker(\inf)\rangle=\langle\Pi\rangle\cap \Ker(\inf)$.
\end{enumerate}
\end{transfer theorem}

More specifically, the intersections $T(G),\bar T(G)$ and the set $\Pi$ are as follows:
We consider a collection of central extensions
\[
\omega:\quad 0\to\dbZ/p\to\dbU_\omega\to\bar\dbU_\omega\to1
\]
of profinite groups subject to certain assumptions (see \S\ref{section on T bar T}).
Let $\alp_\omega$ be the classifying element of $\omega$ in $H^2(\bar\dbU_\omega)$.
For a profinite group $G$ we define $T(G)$ (resp., $\bar T(G)$) to be the intersection of all kernels of continuous homomorphisms $\rho\colon G\to\dbU_\omega$ (resp., $\rho\colon G\to\bar\dbU_\omega$) for some $\omega$.
The set $\Pi$ of special cohomology elements consists of the pullbacks $\bar\rho^*(\alp_\omega)$ of $\alp_\omega$ to $H^2(G/\bar T(G))$ for some $\omega$ and some continuous homomorphism $\bar\rho\colon G/\bar T(G)\to\bar\dbU_\omega$ which lifts to a continuous homomorphism $\rho\colon S\to\dbU_\omega$.
Note that when $S$ is free, the latter condition is always satisfied.
For certain natural choices of $\omega$, these pullbacks yield the cup products, and more generally,  Massey product elements, and Bockstein elements (\S\ref{section on p Zassenhaus filtration}, \S\ref{section on the lower central filtration for odd prime}).
See Theorem \ref{The transfer theorem general form} and Remark \ref{explicit form of the kernel generating condition}.

The paper is organized as follows:
In \S\ref{section on bilinear maps} we assemble some needed facts on bilinear maps.
In \S\ref{section on T bar T}--\S\ref{section on the C-pairing} we develop a formal group-theoretic apparatus which eventually leads to the general Transfer Theorem (Theorem  \ref{The transfer theorem general form}).
It is based on Hoechsmann's theory \cite{Hoechsmann68}, which connects lifting of homomorphisms in group extensions to the transgression map in cohomology.
 This is combined with delicate computations in Pontrjagin duality, which interpret the profinite group side (as in (a)) to the cohomology side (as in (b)).
Then we apply the general theory in three situations:
\begin{enumerate}
\item[(1)]
The \textsl{$p$-Zassenhaus context}, where the central extension arises from $\dbU_n(\dbZ/p)$ and the special cohomology elements are $n$-fold Massey product elements (\S\ref{section on p Zassenhaus filtration});
\item[(2)]
The \textsl{lower $p$-central context}, where the central extensions arise from $\dbU_s(\dbZ/p^{n-s+1})$, $s=1,2\nek n$, and the special cohomology elements are the \textsl{unitriangular spectrum} studied in \cite{Efrat17}, \cite{Efrat20a} (\S\ref{section on lower p-central extension});
\item[(3)]
A \textsl{``mixed" context}, where the central extensions arise from $\dbZ/p^2$ and $M_{p^3}$, $p>2$, and the special cohomology elements are Bockstein elements and sums of Bockstein elements and cup products (\S\ref{section on the lower central filtration for odd prime}).
\end{enumerate}

In each context we obtain a particular cohomological transfer theorem, and moreover, extend earlier results from \cite{Efrat14a}, \cite{CheboluEfratMinac12}, \cite{EfratMinac17}, and \cite{MinacSpira96}.
In \S\ref{section on absolute Galois groups} we apply the results in these three contexts to recover the intersection theorems (i)--(iii) in a uniform manner.
This part uses deep cohomological properties of absolute Galois groups:  The injectivity part of the Merkurjev--Suslin theorem and the description of the Bockstein map as a cup product.
Finally, in \S\ref{section on counterexamples} we give examples which show that two seemingly plausible variants of the intersection theorems fail to hold.
\S\S\ref{section on absolute Galois groups}--\ref{section on counterexamples} are closely related at many points to the papers \cite{EfratMinac11} and \cite{EfratMinac17}, and I thank J\'an Min\'a\v c for the collaboration on these works.

\section{Bilinear maps}
\label{section on bilinear maps}
Let $A$ be a profinite abelian group, $B$ a discrete abelian group, and $Z$ a finite abelian group.
Let $(\cdot,\cdot)\colon A\times B\to Z$ be a bilinear map.
It induces homomorphisms
\[
A\to\Hom(B,Z), \qquad B\to\Hom_{\rm cont}(A,Z).
\]
We refer to their kernels as the \textsl{left  kernel} (resp., the \textsl{right kernel}) of $(\cdot,\cdot)$.
We say that $(\cdot,\cdot)$ is \textsl{left-surjective} (resp., \textsl{right-surjective}) if the left (resp., right) induced homomorphism is surjective.
The map $(\cdot,\cdot)$ is \textsl{non-degenerate} if both induced homomorphisms are injective, and is \textsl{perfect} if they are both isomorphisms.

Let $\alp\colon A_1\to A_2$ be a continuous homomorphism of profinite abelian groups and $\beta\colon B_2\to B_1$ a homomorphism of discrete abelian groups.
Suppose that we have a diagram of bilinear maps
\begin{equation}
\label{cd of bilinear maps}
\xymatrix{
A_1\ar[d]_{\alp}&*-<3pc>{\times}&B_1\ar[r]^{(\cdot,\cdot)_1}&Z\ar@{=}[d]\\
A_2&*-<3pc>{\times}&B_2\ar[u]_{\beta}\ar[r]^{(\cdot,\cdot)_2}&Z
}
\end{equation}
which is commutative, in the sense that $(a,\beta(b))_1=(\alp(a),b)_2$ for every $a\in A_1$ and $b\in B_2$.

\begin{lem}
\label{Coker Ker}
\begin{enumerate}
\item[(a)]
The commutative diagram (\ref{cd of bilinear maps}) induces a bilinear map $(\cdot,\cdot)_{\Coker,\Ker}\colon \Coker(\alp)\times\Ker(\beta)\to Z$.
\item[(b)]
If the right kernel of $(\cdot,\cdot)_2$ is trivial, then the right kernel of $(\cdot,\cdot)_{\Coker,\Ker}$ is trivial.
\item[(c)]
If the natural map $\pi\colon A_2\to\Coker(\alp)$ splits, the right kernel of $(\cdot,\cdot)_1$ is trivial,  and $(\cdot,\cdot)_2$ is right-surjective, then $(\cdot,\cdot)_{\Coker,\Ker}$ is right-surjective.
\end{enumerate}
\end{lem}
\begin{proof}
For $a'\in A_2$ and $b\in\Ker(\beta)$ we set $(\pi(a'),b)_{\Coker,\Ker}=(a',b)_2$ (with $\pi$ as in (c)).
It is well-defined, since for $a'=\alp(a)$, $a\in A_1$, we have $(a',b)_2=(a,\beta(b))_1=(a,0)_1=0$.
It is straightforward to verify the bilinearity, as well as assertion (b).

For (c), we decompose $A_2=\Img(\alp)\oplus A'_2$, where $\pi$ maps $A'_2$ bijectively onto $\Coker(\alp)$.
Let $\lam_0\colon \Coker(\alp)\to Z$ be a continuous homomorphism.
Let $\lam\colon A_2\to Z$ be the continuous homomorphism which is $0$ on $\Img(\alp)$ and is $\lam_0\circ\pi$ on $A'_2$.
By assumption, there exists $b\in B_2$ such that $\lam=(\cdot, b)_2$.
For every $a\in A_1$ we have $(a,\beta(b))_1=(\alp(a),b)_2=\lam(\alp(a))=0$, whence $b\in\Ker(\beta)$.
Moreover, for every $a'\in A'_2$ one has
\[
\lam_0(\pi(a'))=\lam(a')=(a',b)_2=(\pi(a'),b)_{\Coker,\Ker},
\]
i.e., $\lam_0=(\cdot,b)_{\Coker,\Ker}$, as required.
\end{proof}

\begin{cor}
\label{Coker Ker is perfect}
Suppose that $A_i,B_i$, $i=1,2$, have prime exponent $p$, and that $Z=\dbZ/p$.
Assume that the right kernel of $(\cdot,\cdot)_1$ is trivial and that $(\cdot,\cdot)_2$ is perfect.
Then $(\cdot,\cdot)_{\Coker,\Ker}$ is perfect.
\end{cor}
\begin{proof}
Since $A_2$ is an elementary abelian $p$-group, the splitting assumption  in Lemma \ref{Coker Ker}(c) is satisfied.
Hence  $(\cdot,\cdot)_{\Coker,\Ker}$ is right-surjective.
By Lemma \ref{Coker Ker}(b), its right kernel is trivial, so the induced map $\Ker(\beta)\to\Hom_{\rm cont}(\Coker(\alp),\dbZ/p)$ is an isomorphism.
Thus $\Coker(\alp)$ and $\Ker(\beta)$ are Pontrjagin duals, and the assertion follows.
\end{proof}

\begin{lem}
\label{left-surjectivity}
Suppose that $A_i,B_i$, $i=1,2$, have prime exponent $p$, and that $Z=\dbZ/p$.
Assume that $(\cdot,\cdot)_1$ is left-surjective and $\beta$ is injective.
Then $(\cdot,\cdot)_2$ is left-surjective.
\end{lem}
\begin{proof}
The assumptions imply that $B_2$ is a direct factor of $B_1$.
Hence the dual map $\beta^\vee\colon \Hom(B_1,\dbZ/p)\to \Hom(B_2,\dbZ/p)$ is surjective.
Therefore in  the commutative square
\[
\xymatrix{
A_1\ar@{->>}[r]\ar[d]_{\alpha}&\Hom(B_1,\dbZ/p)\ar@{->>}[d]^{\beta^\vee}\\
A_2\ar[r]&\Hom(B_2,\dbZ/p).
}
\]
the lower horizontal map is surjective.
\end{proof}

\section{The subgroups $T$ and $\bar T$}
\label{section on T bar T}
For profinite groups $\dbU$ and $G$, let
\[
T^\dbU(G)=\bigcap\Ker(\rho),
 \]
 where the intersection is over all continuous homomorphisms $\rho\colon G\to\dbU$.
Alternatively, $T^\dbU(G)=\bigcap N$, where $N$ ranges over all closed normal subgroups of $G$ such that $G/N$ embeds as a closed subgroup of $\dbU$.
The proof of the following lemma is straightforward:

\begin{lem}
\label{intersection of kernels}
Let $\dbU,\dbU_1\nek\dbU_n$ be profinite groups, and suppose that there are continuous homomorphisms $\gam_i\colon\dbU\to\dbU_i$, $i=1,2\nek n$, such that $\bigcap_{i=1}^n\Ker(\gam_i)=\{1\}$.
For every profinite group $G$ one has $\bigcap_{i=1}^nT^{\dbU_i}(G)\leq T^\dbU(G)$.
\end{lem}

Let $Z$ be a finite abelian group.
We fix a nonempty set $\Omega$ of central extensions of profinite groups
\[
\omega:\quad 0\to Z\to\dbU\xrightarrow{\lam}\bar\dbU\to1.
\]
We assume that for every $\omega$:
\begin{enumerate}
\item[(I)]
There is a family of continuous homomorphisms $\gam_i\colon\bar\dbU\to\dbU$, $i=1,2\nek n$, such that $\bigcap_{i=1}^n\Ker(\gam_i)=\{1\}$.
\item[(II)]
$Z$ embed in $\bar\dbU$.
\end{enumerate}
For a profinite group $G$, assumption (I) and Lemma \ref{intersection of kernels} imply that
$T^\dbU(G)\leq T^{\bar\dbU}(G)$.
Let
\[
T(G)=\bigcap_{\omega\in\Omega}T^{\dbU_\omega}(G),\quad \bar T(G)=\bigcap_{\omega\in\Omega}T^{\bar\dbU_\omega}(G).
\]
Then $T(G)$, $\bar T(G)$ are closed normal subgroups of $G$, and
\[
T(G)\leq \bar T(G).
\]

If $\pi\colon G_1\to G_2$ is a continuous homomorphism of profinite groups, then
\begin{equation}
\label{weak functoriality of T, bar T}
\pi[T(G_1)]\leq T(G_2), \qquad \pi[\bar T(G_1)]\leq \bar T(G_2).
\end{equation}

\begin{lem}
\label{functoriality of T, bar T}
Let $\pi\colon G_1\to G_2$ is a continuous epimorphism of profinite groups.
Then:
\begin{enumerate}
\item[(a)]
$\Ker(\pi)\leq T(G_1)$ if and only if $T(G_1)=\pi\inv[T(G_2)]$.
\item[(b)]
$\Ker(\pi)\leq \bar T(G_1)$ if and only if $\bar T(G_1)=\pi\inv[\bar T(G_2)]$.
\end{enumerate}
\end{lem}
\begin{proof}
(a) \quad
The ``if" part is immediate.

For the ``only if" part, assume that $\Ker(\pi)\leq T(G_1)$.
Then every continuous homomorphism $\rho_1\colon G_1\to\dbU_\omega$, $\omega\in\Omega$, factors via a continuous homomorphism $\rho_2\colon G_2\to\dbU_\omega$.
We have $\rho_1[\pi\inv[T(G_2)]]=\rho_2[T(G_2)]=1$.
Thus $\pi\inv[T(G_2)]\leq T(G_1)$.
The opposite inclusion follows from (\ref{weak functoriality of T, bar T}).

\medskip

(b)\quad
Similarly.
\end{proof}

\begin{lem}
\label{exponent}
If $Z$ has exponent $m$, then $\bar T(G)^m[G,\bar T(G)]\leq  T(G)$ for every profinite group $G$.
Consequently, $\bar T(G)/T(G)$ is an abelian group of exponent dividing $m$.
\end{lem}
\begin{proof}
Let $\omega\in\Omega$.
Since it is central, $Z^m[\dbU_\omega,Z]=1$.
For every continuous homomorphism $\rho\colon G\to \dbU_\omega$ we have $\rho(\bar T(G))\leq \Ker(\lam_\omega)=Z$, whence $\rho(\bar T(G)^m[G,\bar T(G)])=1$, and the assertion follows.
\end{proof}

For a profinite group $G$ which acts trivially on $Z$, we write $H^i(G)=H^i(G,Z)$ for the $i$th (profinite) cohomology group.
Recall that $H^1(G)=\Hom_{\rm cont}(G,Z)$.
We write $\inf$, $\res$ and $\trg$ for the inflation, restriction, and transgression homomorphisms, respectively \cite{NeukirchSchmidtWingberg}*{Ch.\ I}.

\begin{lem}
\label{inflation is an isomorphism}
Let $N$ be a closed normal subgroup of $G$ contained in $\bar T(G)$.
\begin{enumerate}
\item[(a)]
The map $\inf\colon H^1(G/N)\to H^1(G)$ is an isomorphism.
\item[(b)]
There is an exact sequence
\[
0\to H^1(N)^G\xrightarrow{\trg} H^2(G/N)\xrightarrow{\inf}H^2(G).
\]
\end{enumerate}
\end{lem}
\begin{proof}
(a) \quad
Let $\rho\colon G\to Z$ be a continuous homomorphism.
By assumption (II), $Z$ embeds in every $\bar\dbU_\omega$, so $\rho(\bar T(G))=1$, whence $\rho(N)=1$.
Therefore $\rho$ factors via a unique continuous homomorphism $G/N\to Z$.

\medskip

(b) \quad
This follows from (a) and the five term sequence in profinite cohomology  \cite{NeukirchSchmidtWingberg}*{Prop.\ 1.6.7}.
\end{proof}

\begin{rems}
\label{functoriality of 5-term sequence}
\rm
(1) \quad
The exact sequence in Lemma \ref{inflation is an isomorphism}(b) is functorial.
Namely, for every continuous homomorphism $\pi\colon G\to G'$ of profinite groups, and every closed normal subgroups $N$ of $G$ and $N'$ of $G'$ such that $N\leq \bar T(G)$, $N'\leq\bar T(G')$, and $\pi[N]\leq N'$, there is a commutative diagram
with exact rows
\begin{equation}
\label{functoriality of 5 term sequence}
\xymatrix{
0\ar[r]&H^1(N')^{G'}\ar[r]^{\trg}\ar[d]& H^2(G'/N')\ar[r]^{\quad\inf}\ar[d] &H^2(G')\ar[d]\\
0\ar[r]&H^1(N)^G\ar[r]^{\trg}& H^2(G/N)\ar[r]^{\quad\inf}&H^2(G),
}
\end{equation}
where the vertical maps are induced by $\pi$.

\medskip

(2) \quad
In particular, when $G=G'$ and $\pi=\id$, the snake lemma yields an isomorphism between  the kernels of the left and middle vertical maps in (\ref{functoriality of 5 term sequence}).
\end{rems}

\section{Liftable homomorphisms}
\label{section on liftable homomorphisms}
Let $G$ be a profinite group, let $N$ be a closed normal subgroup of $G$ contained in $\bar T(G)$, and let $\pi\colon G\to G/N$ be the natural epimorphism.
Consider $\omega\in\Omega$ and a continuous homomorphism $\bar\rho\colon G/N\to\bar\dbU_\omega$.
We say that $\bar\rho$ is a \textsl{$\pi$-liftable} homomorphism if there exists a continuous homomorphism $\rho\colon G\to\dbU_\omega$ with a commutative square
\[
\xymatrix{
G\ar[r]^{\pi\ }\ar[d]_{\rho}&G/N\ar[d]^{\bar\rho}\\
\dbU_\omega\ar[r]^{\lam_\omega}&\ \bar\dbU_\omega.
}
\]
We write $\dbU_\omega\times_{\bar\dbU_\omega}(G/N)$ for the fiber product with respect to $\lam_\omega$ and $\bar\rho$.

Let $\alp_\omega\in H^2(\bar\dbU_\omega)$ be the classifying element of the extension $\omega$ under the Schreier correspondence \cite{NeukirchSchmidtWingberg}*{Th.\ 1.2.4}.
Let $\bar\rho^*(\alp_\omega)$ be the pullback of $\alp_\omega$ to $H^2(G/N)$ via $\bar\rho$.
It corresponds to the central extension
\[
0\to Z\to\dbU_\omega\times_{\bar\dbU_\omega}(G/N)\to G/N\to1,
\]
where the right map is the projection map \cite{GilleSzamuely06}*{Remark 3.3.11}.

Now consider the following commutative diagram:
\begin{equation}
\label{diagram for liftable homs}
\xymatrix{
&&&G\ar[ddl]^(0.7){\rho}\ar[d]^{\pi}\ar[dl]_{\Psi}&\\
0\ar[r]&Z\ar[r]\ar@{=}[d]&\dbU_\omega\times_{\bar\dbU_\omega}(G/N)\ar[r]\ar[d]&G/N\ar[d]^{\bar\rho}\ar[r]&1\\
0\ar[r]&Z\ar[r]&\dbU_\omega\ar[r]^{\lam_\omega}&\bar\dbU_\omega\ar[r]&1.
}
\end{equation}

\begin{prop}
\label{equivalences for liftable homs}
In this set-up, the following conditions are equivalent:
\begin{enumerate}
\item[(a)]
$\bar\rho$ is $\pi$-liftable, i.e., there exists a continuous homomorphism $\rho\colon G\to \dbU_\omega$ making the lower triangle commutative;
\item[(b)]
There exists a continuous homomorphism $\Psi\colon G\to \dbU_\omega\times_{\bar\dbU_\omega}(G/N)$ making the upper triangle commutative;
\item[(c)]
The inflation of $\bar\rho^*(\alp_\omega)$ to $H^2(G)$ is trivial;
\item[(d)]
There is a (necessarily unique)  $\psi\in H^1(N)^G$ with $\bar\rho^*(\alp_\omega)=\trg(\psi)$.
\end{enumerate}
\end{prop}
\begin{proof}
(a)$\Leftrightarrow$(b): \quad
Use the universal property of the fiber product.

\medskip

(b)$\Leftrightarrow$(c):  \quad
This follows from Hoechsmann's lemma (\cite{Hoechsmann68}*{1.1}, \cite{NeukirchSchmidtWingberg}*{Prop.\ 3.5.9}).

\medskip

(c)$\Leftrightarrow$(d): \quad
Apply Lemma \ref{inflation is an isomorphism}(b).
\end{proof}

\begin{rems}
\label{remarks on liftable}
\rm
(1) \quad
In this situation, $\Ker(\Psi)=\Ker(\rho)\cap N$.

\medskip

(2) \quad
Since $N\leq\bar T(G)$, every continuous homomorphism $\rho\colon G\to\dbU_\omega$, with $\omega\in \Omega$, induces a $\pi$-liftable $\bar\rho\colon G/N\to\bar\dbU_\omega$.

\medskip

(3) \quad
If $\cd_p(G)\leq1$ for every prime divisor $p$ of the order of $Z$, then the embedding problem in the lower triangle of (\ref{diagram for liftable homs}) is solvable \cite{NeukirchSchmidtWingberg}*{Th.\ 3.5.6}, i.e., $\bar\rho$ is $\pi$-liftable.

\medskip
(4) \quad
By the explicit description of $\trg$ \cite{NeukirchSchmidtWingberg}*{Prop.\ 1.6.6}, the homomorphism $\psi$ in (d) is the restriction $\Psi|_N$.
Considering it as a map to the left factor of $\dbU_\omega\times_{\bar\dbU_\omega}(G/N)$, we thus have $\psi=\rho|_N$.

\medskip

(5) \quad
Every continuous homomorphism $\rho\colon G\to\dbU_\omega$ as above factors via $G/(N\cap T(G))$.
Hence $\bar\rho$ is $\pi$-liftable if and only if it is $\pi'$-liftable, where $\pi'\colon G/(N\cap T(G))\to G/N$ is the natural map.
Consequently, by the equivalence (a)$\Leftrightarrow$(c) in Proposition \ref{equivalences for liftable homs}, $\bar\rho$ is $\pi$-liftable if and only if $\inf(\bar\rho^*(\alp_\omega))=0$ in $H^2(G/(N\cap T(G))$.
\end{rems}

\begin{defin}
\label{def of liftable part of cohomology}
\rm
We call an element $\bar\rho^*(\alp_\omega)$ of $H^2(G/N)$, where $\omega\in\Omega$ and $\bar\rho\colon G/N\to\bar\dbU_\omega$ is a $\pi$-liftable continuous homomorphism, a \textsl{$\pi$-liftable pullback}.

The \textsl{$\pi$-liftable part} of the cohomology group $H^2(G/N)$ is its subgroup $H^2(G/N)_\pi$ generated by all $\pi$-liftable pullbacks $\bar\rho^*(\alp_\omega)$.
\end{defin}

\section{The $A$-pairing}
From now on we assume  that
\[
Z\isom \dbZ/m
\]
for an integer $m\geq2$.
Thus for a profinite group $G$ we set $H^i(G)=H^i(G,\dbZ/m)$.
For a closed normal subgroup $N$ of $G$, the conjugation induces a $G$-action on $H^1(N)$ \cite{NeukirchSchmidtWingberg}*{Ch.\ I, \S5}.
The $G$-invariant part of $H^1(N)$ satisfies
\[
H^1(N)^G\isom\Hom_{\rm cont}(N/N^m[G,N],\dbZ/m),
\]
which is the Pontrjagin dual $(N/N^m[G,N])^\vee$ of the torsion abelian group $N/N^m[G,N]$.
Thus the substitution map
\[
N/N^m[G,N]\times H^1(N)^G\to\dbZ/m
\]
is a perfect bilinear map.

For closed normal subgroups $N_1\leq N_2$ of $G$ there is a commutative diagram of perfect bilinear maps
\[
\xymatrix{
N_1/N_1^m[G,N_1]\ar[d]& *-<3pc>{\times}& H^1(N_1)^G\ar[r]&\dbZ/m\ar@{=}[d]\\
N_2/N_2^m[G,N_2]&*-<3pc>{\times}& H^1(N_2)^G\ar[u]^{\res}\ar[r]&\dbZ/m.\\
}
\]
By Lemma \ref{Coker Ker}(b), it gives rise to a bilinear map with trivial right kernel
\begin{equation}
\label{cokernel of substitution pairings}
N_2/N_1N_2^m[G,N_2]\times\Ker\bigl(H^1(N_2)^G\xrightarrow{\res}H^1(N_1)\bigr)\to \dbZ/m.
\end{equation}
We set
\[
A_G(N_1,N_2)=\Ker\bigl(H^2(G/N_2)\xrightarrow{\inf} H^2(G/N_1)\bigr).
\]
By Remark \ref{functoriality of 5-term sequence}(2), these two kernels are isomorphic via transgression, so (\ref{cokernel of substitution pairings}) induces a bilinear map with trivial right kernel
\[
\langle\cdot,\cdot\rangle_{N_1,N_2}^A\colon  N_2/N_1N_2^m[G,N_2]\times A_G(N_1,N_2)\to \dbZ/m.
\]
Moreover, Corollary \ref{Coker Ker is perfect} gives:

\begin{prop}
\label{transgression perfect pairing}
When $m=p$ is prime, $\langle\cdot,\cdot\rangle_{N_1,N_2}^A$ is perfect.
\end{prop}

\begin{rem}
\label{functoriality of the transgression pairing}
\rm
This construction is functorial.
Namely, let $\pi\colon G\to G'$ be a continuous homomorphism of profinite groups.
Let $N_1\leq N_2$ be closed normal subgroups of $G$, let $N'_1\leq N'_2$ be closed normal subgroups of $G'$, and suppose that $\pi[N_1]\leq N'_1$ and $\pi[N_2]\leq N'_2$.
Then $\pi$ induces a commutative diagram
\[
\xymatrix{
N_2/N_1N_2^m[G,N_2]\ar[d]&*-<3pc>{\times}& A_G(N_1,N_2)\ar[rr]^{\quad\langle\cdot,\cdot\rangle_{N_1,N_2}^A}&&\dbZ/m\ar@{=}[d]\\
N'_2/N'_1(N'_2)^m[G',N'_2]&*-<3pc>{\times}& A_{G'}(N'_1,N'_2)\ar[u]\ar[rr]^{\quad\langle\cdot,\cdot\rangle_{N'_1,N'_2}^A}&&\dbZ/m.\\
}
\]
\end{rem}

\section{The $B$-pairing}
\label{section on pairings for B}
Let $G$ be a profinite group, and let $N_1\leq N_2$ be closed normal subgroups of $G$ contained in $\bar T(G)$.
Let $\pi_i\colon G\to G/N_i$, $i=1,2$, be the natural projections.
Every $\pi_1$-liftable continuous homomorphism $\bar\rho\colon G/N_1\to\bar\dbU_\omega$, $\omega\in\Omega$, factors via a $\pi_2$-liftable continuous homomorphism $\bar\rho\colon G/N_2\to\bar\dbU_\omega$.
Hence the inflation homomorphism maps $H^2(G/N_2)_{\pi_2}$ \textsl{onto} $H^2(G/N_1)_{\pi_1}$ (see Definition \ref{def of liftable part of cohomology}).
We set
\[
B_G(N_1,N_2)=\Bigl\langle \bar\rho^*(\alp_\omega)\in A_G(N_1,N_2) \  \Bigm|  \  \omega\in\Omega,\ \bar\rho\colon G/N_2\to\bar\dbU_\omega  \hbox{  $\pi_2$-liftable}\Bigr\rangle.
\]
Thus $B_G(N_1,N_2)$ is generated by all $\pi_2$-liftable pullbacks in $ H^2(G/N_2)_{\pi_2}$ which vanish under the inflation to $H^2(G/N_1)$.

\begin{thm}
\label{pairing for B}
In this setup, there is a non-degenerate bilinear map
\[
\langle\cdot,\cdot\rangle_{N_1,N_2}^B\colon N_2/(N_2\cap \pi_1\inv[T(G/N_1)])\times B_G(N_1,N_2)\to\dbZ/m
\]
which fits into a commutative diagram
\[
\xymatrix{
N_2/N_1N_2^m[G,N_2]\ar@{->>}[d]&*-<3pc>{\times}&A_G(N_1,N_2)\ar[rr]^{\quad\langle\cdot,\cdot\rangle_{N_1,N_2}^A}&&\dbZ/m\ar@{=}[d]\\
N_2/(N_2\cap \pi_1\inv[T(G/N_1)])&*-<3pc>{\times}&B_G(N_1,N_2)\ar@{^{(}->}[u]\ar[rr]^{\quad\langle\cdot,\cdot\rangle_{N_1,N_2}^B}&&\dbZ/m.
}
\]
Moreover, when $m=p$ is prime, $\langle\cdot,\cdot\rangle_{N_1,N_2}^B$ is perfect.
\end{thm}
\begin{proof}
First we note that, by  Lemma \ref{exponent} and (\ref{weak functoriality of T, bar T}),
\begin{equation}
\label{inclusion of subgroups}
\begin{split}
N_1N_2^m[G,N_2]\leq N_2\cap N_1\bar T(G)^m[G,\bar T(G)]&\leq N_2\cap N_1T(G)\\
&\leq N_2\cap \pi_1\inv[T(G/N_1)].
\end{split}
\end{equation}

Now let $B_0$ be the above set of generators of $B_G(N_1,N_2)$, i.e., $B_0$ is the set of all $\pi_2$-liftable pullbacks $\bar\rho^*(\alp_\omega)$ such that
$\inf_{G/N_1}(\bar\rho^*(\alp_\omega))=0$.
By Proposition \ref{equivalences for liftable homs}, the latter condition means that $\bar\rho^*(\alp_\omega)$ is in fact $\pi_{12}$-liftable, where $\pi_{12}\colon G/N_1\to G/N_2$ is the natural projection.
Equivalently, there is a (necessarily unique) $\psi\in H^1(N_2/N_1)^{G/N_1}$ such that $\bar\rho^*(\alp_\omega)=\trg(\psi)$, for the transgression map $\trg\colon  H^1(N_2/N_1)^{G/N_1}\to H^2(G/N_2)$.
For a $\pi_{12}$-lifting $\rho$ of $\bar\rho$, we have $\rho|_{N_2/N_1}=\psi$ (Remark \ref{remarks on liftable}(4)).

Next let $\sig\in N_2$, and let  $\bar\sig$ be its coset modulo $N_1N_2^m[G,N_2]$.
Then $\bar\sig$ is in the annihilator of $B_G(N_1,N_2)$ with respect to $\langle\cdot,\cdot\rangle_{N_1,N_2}^A$ if and only if it is annihilated by all generators in $B_0$.
This means that for every $\bar\rho,\psi,\rho$ as in the previous paragraph,
\[
0=\langle\bar\sig,\bar\rho^*(\alp_\omega)\rangle_{N_1,N_2}^A=\langle\bar\sig,\trg(\psi)\rangle_{N_1,N_2}^A=\psi(\sig N_1)=\rho(\sig N_1).
\]
Equivalently,  $\sig\in N_2\cap\pi_1\inv[T(G/N_1)]$.

Consequently, $\langle\cdot,\cdot\rangle_{N_1,N_2}^A$ induces a bilinear map $\langle\cdot,\cdot\rangle_{N_1,N_2}^B$ as in the assertion, whose left kernel is trivial.
In light of (\ref{inclusion of subgroups}), this bilinear map fits into a commutative diagram as claimed.
Since the right kernel of $\langle\cdot,\cdot\rangle_{N_1,N_2}^A$ is trivial and the left (resp., right) vertical map of the diagram is surjective (resp., injective), the right kernel of $\langle\cdot,\cdot\rangle_{N_1,N_2}^B$ is also trivial, i.e., $\langle\cdot,\cdot\rangle_{N_1,N_2}^B$ is non-degenerate.

When $m=p$ is prime, $\langle\cdot,\rangle_{N_1,N_2}^A$ is perfect, whence left-surjective.
It follows from Lemma \ref{left-surjectivity} that $\langle\cdot,\cdot\rangle_{N_1,N_2}^B$ is also left-surjective.
Since it is non-degenerate, Pontrjagin duality implies that it is perfect.
\end{proof}

\begin{cor}
\label{an exact sequence}
In the setup of Theorem \ref{pairing for B}, and for $m=p$ prime, the following conditions are equivalent:
\begin{enumerate}
\item[(a)]
$N_1N_2^p[G,N_2]=N_2\cap\pi_1\inv[T(G/N_1)]$;
\item[(b)]
There is an exact sequence
\[
0\to B_G(N_1,N_2)\hookrightarrow H^2(G/N_2)\xrightarrow{\inf}H^2(G/N_1).
\]
\end{enumerate}
\end{cor}

\begin{rem}
\label{functoriality of the liftable pairing}
\rm
The bilinear map $\langle\cdot,\cdot\rangle_{N_1,N_2}^B$ is functorial in the following sense:

Let $\pi\colon G\to G'$ be a continuous map of profinite groups.
Let $N_1\leq N_2$ be closed normal subgroups of $G$ contained in $\bar T(G)$, let $N'_1\leq N'_2$ be closed normal subgroups of $G'$ contained in $\bar T(G')$, and suppose that $\pi[N_1]\leq N'_1$ and $\pi[N_2]\leq N'_2$.
We write $\pi_i\colon G\to G/N_i$, $\pi'_i\colon G'\to G'/N'_i$, $i=1,2$, for the natural projections.
Then the functoriality of $\langle\cdot,\cdot\rangle_{N_1,N_2}^A$ and $\langle\cdot,\cdot\rangle_{N'_1,N'_2}^A$ (Remark \ref{functoriality of the transgression pairing}) and the construction of the pairing in Theorem \ref{pairing for B} give rise to a commutative diagram of non-degenerate bilinear maps
\[
\xymatrix{
N_2/(N_2\cap\pi_1\inv[T(G/N_1)]\ar[d]&*-<3pc>{\times}& B_G(N_1,N_2)\ar[rr]^{\ \langle\cdot,\cdot\rangle_{N_1,N_2}^B}&&\dbZ/m\ar@{=}[d]\\
N'_2/(N'_2\cap(\pi'_1)\inv[T(G'/N'_1)]&*-<3pc>{\times}& B_{G'}(N'_1,N'_2)\ar[u]\ar[rr]^{\ \langle\cdot,\cdot\rangle_{N'_1,N'_2}^B}&&\dbZ/m.
}
\]
\end{rem}

As a special case we obtain:

\begin{prop}
\label{a special case}
Let $N$ be a closed normal subgroup of $G$ contained in $\bar T(G)$, and let $\pi\colon G\to G/N$ be the natural projection.
\begin{enumerate}
\item[(a)]
There is a non-degenerate bilinear map
\[
N/(N\cap T(G))\times H^2(G/N)_\pi\to\dbZ/m.
\]
\item[(b)]
When $m=p$ is prime, there is an exact sequence
\[
0\to H^2(G/N)_\pi\to H^2(G/N)\xrightarrow{\inf}H^2(G/(N\cap T(G))).
\]
\end{enumerate}
\end{prop}
\begin{proof}
Take in Theorem \ref{pairing for B} $N_1=N\cap T(G)$ and $N_2=N$.
By Remark \ref{remarks on liftable}(5), $\inf\colon H^2(G/N_2)_\pi\to H^2(G/N_1)$ is trivial, so $B_G(N_1,N_2)= H^2(G/N)_\pi$.
Further, by Lemma \ref{functoriality of T, bar T},  $\pi_1\inv[T(G/N_1)]=T(G)$, and we deduce (a).

By Lemma \ref{exponent}, $N^m[G,N]\leq N\cap T(G)=N_1$, so $N_1N^m[G,N]=N\cap T(G)$.
Therefore Corollary \ref{an exact sequence} gives (b).
\end{proof}

\section{The $C$-pairing}
\label{section on the C-pairing}
Let $N_1\leq N_2$ be again closed normal subgroups of the profinite group $G$ which are contained in $\bar T(G)$.
Let $\pi_i\colon G\to G/N_1$, $i=1,2$, be the canonical projections.
We set
\[
C_G(N_1,N_2)=\Ker\bigl(\inf\colon H^2(G/N_2)_{\pi_2}\to H^2(G/N_1)_{\pi_1}\bigr).
\]
Then
\[
B_G(N_1,N_2)\leq C_G(N_1,N_2)\leq A_G(N_1,N_2).
\]

\begin{thm}
\label{pairing for C}
Suppose that $m=p$ is prime.
Then $\langle\cdot,\cdot\rangle^B_{N_1,N_2}$ induces a perfect bilinear map
\[
\langle\cdot,\cdot\rangle_{N_1,N_2}^C\colon  N_2/(N_2\cap N_1T(G))\times  C_G(N_1,N_2)\to\dbZ/p.
\]
\end{thm}
\begin{proof}
Let
\[
\widetilde\pi_1\colon G\to G/(N_1\cap T(G)), \quad  \widetilde\pi_2\colon G\to G/(N_2\cap T(G))
\]
be the canonical projections.
The functoriality of $\langle\cdot,\cdot\rangle^B$ (Remark \ref{functoriality of the liftable pairing}) gives rise to a commutative diagram of perfect bilinear maps
\[
\xymatrix{
N_1/(N_1\cap \widetilde\pi_1\inv[T(G/(N_1\cap T(G))]) \ar[d]&*-<3pc>{\times}& B_G(N_1\cap T(G),N_1)\ar[r]^{\quad\qquad\langle\cdot,\cdot\rangle^B}&\dbZ/p\ar@{=}[d]\\
N_2/(N_2\cap \widetilde\pi_2\inv[T(G/(N_2\cap T(G))] )&*-<3pc>{\times}& B_G(N_2\cap T(G),N_2)\ar[u]\ar[r]^{\quad\qquad\langle\cdot,\cdot\rangle^B}&\dbZ/p,
}
\]
(where we omit the obvious subscripts on $\langle\cdot,\cdot\rangle^B$).
By Lemma \ref{functoriality of T, bar T},
\[
\widetilde\pi_1\inv[T(G/(N_1\cap T(G)))]=T(G), \quad \widetilde\pi_2\inv[T(G/(N_2\cap T(G)))]=T(G).
\]
In view of  Remark \ref{remarks on liftable}(5),
\[
B_G(N_1\cap T(G),N_1)=H^2(G/N_1)_{\pi_1}, \quad B_G(N_2\cap T(G),N_2)=H^2(G/N_2)_{\pi_2}.
\]
We obtain a commutative diagram of perfect bilinear maps
\[
\xymatrix{
N_1/(N_1\cap T(G)) \ar[d]&*-<3pc>{\times}& H^2(G/N_1)_{\pi_1}\ar[r]^{\qquad\langle\cdot,\cdot\rangle^B}&\dbZ/p\ar@{=}[d]\\
N_2/(N_2\cap T(G))&*-<3pc>{\times}& H^2(G/N_2)_{\pi_2}\ar[u]\ar[r]^{\qquad\langle\cdot,\cdot\rangle^B}&\dbZ/p.
}
\]
By Corollary \ref{Coker Ker is perfect}, it induces a perfect bilinear map between the cokernel $N_2/N_1(N_2\cap T(G))\isom N_2/(N_2\cap N_1T(G))$ of the left vertical map and the kernel $C_G(N_1,N_2)$ of the right vertical map, as required.
\end{proof}

\begin{rem}
\rm
The bilinear map $\langle\cdot,\cdot\rangle^C_{N_1,N_2}$ is functorial in a similar manner to $\langle\cdot,\cdot\rangle^A_{N_1,N_2}$, $\langle\cdot,\cdot\rangle^B_{N_1,N_2}$  (see Remarks \ref{functoriality of the transgression pairing} and \ref{functoriality of the liftable pairing}).
\end{rem}

\begin{cor}
\label{generalization of a result of Dwyer}
The following conditions are equivalent:
\begin{enumerate}
\item[(a)]
$\inf\colon H^2(G/N_2)_{\pi_2}\to H^2(G/N_1)_{\pi_1}$ is an isomorphism;
\item[(b)]
$\inf\colon H^2(G/N_2)_{\pi_2}\to H^2(G/N_1)_{\pi_1}$ is a monomorphism;
\item[(c)]
$N_1T(G)=N_2T(G)$.
\end{enumerate}
\end{cor}
\begin{proof}
We have noticed in \S\ref{section on pairings for B} that $\inf$ is surjective, so (a) and (b) are equivalent.
The equivalence of (b) and (c) follows from Theorem \ref{pairing for C}.
\end{proof}

We summarize Theorem \ref{pairing for B} (for $m=p$ prime) and Theorem \ref{pairing for C}  by the following commutative diagram of perfect bilinear maps.
Here $N_1\leq N_2$ are closed normal subgroups of $G$ contained in $\bar T(G)$:
\begin{equation}
\label{summary}
\xymatrix{
N_2/N_1N_2^p[G,N_2]\ar@{->>}[d]&*-<3pc>{\times}& A_G(N_1,N_2)\ar[rr]^{\quad\langle\cdot,\cdot\rangle_{N_1,N_2}^A}&&\dbZ/p\ar@{=}[d]\\
N_2/(N_2\cap N_1T(G))\ar@{->>}[d]&*-<3pc>{\times}& C_G(N_1,N_2)\ar@{^{(}->}[u]\ar[rr]^{\quad\langle\cdot,\cdot\rangle_{N_1,N_2}^C}&&\dbZ/p\ar@{=}[d]\\
N_2/(N_2\cap\pi_1\inv[T(G/N_1)])&*-<3pc>{\times}& B_G(N_1,N_2)\ar@{^{(}->}[u]\ar[rr]^{\quad\langle\cdot,\cdot\rangle_{N_1,N_2}^B}&&\dbZ/p.
}
\end{equation}

\begin{defin}
\rm
We say that the \textsl{kernel generating condition} holds for the subgroups $N_1,N_2$ of $G$ if $B_G(N_1,N_2)=C_G(N_1,N_2)$.
\end{defin}

We now arrive at the main result of the construction, which characterizes the transfer condition $\pi[T(G)]=T(\pi[G])$ in cohomological terms:

\begin{thm}
\label{The transfer theorem general form}
Let $m=p$ be prime.
Let $N$ be a closed normal subgroup of the profinite group $G$ contained in $\bar T(G)$, and let $\pi\colon G\to G/N$ be the natural projection.
The following conditions are equivalent:
\begin{enumerate}
\item[(a)]
$\pi[T(G)]=T(G/N)$;
\item[(b)]
The kernel generating condition holds for the subgroups $\bar T(G)$, $N$ of $G$.
\end{enumerate}
\end{thm}
\begin{rem}
\label{explicit form of the kernel generating condition}
\rm
Condition (b) can be stated more explicitly as follows:
Let $\bar\pi\colon G\to G/\bar T(G)$ be the natural projection, and consider the inflation map $\inf\colon H^2(G/\bar T(G))\to H^2(G/N)$.
Let $\Pi$ be the set of all pullbacks $\bar\rho^*(\alp_\omega)$, with $\omega\in\Omega$ and $\bar\rho\colon  G/\bar T(G)\to \bar\dbU_\omega$ is a $\bar\pi$-liftable continuous homomorphism.
Then the subgroups $\langle\Pi\rangle$, $\langle\Pi\cap\Ker(\inf)\rangle$ of $H^2(G/\bar T(G))$ generated by $\Pi$, $\Pi\cap\Ker(\inf)$, respectively, satisfy
\[
\langle\Pi\cap\Ker(\inf)\rangle=\langle\Pi\rangle\cap\Ker(\inf).
\]
\end{rem}

\begin{proof}
First note that (a) is equivalent to $NT(G)=\pi\inv[T(G/N)]$.
As $N\leq\bar T(G)$, Lemma \ref{functoriality of T, bar T} implies that
\[
\pi\inv[T(G/N)]\leq\pi\inv[\bar T(G/N)]=\bar T(G).
\]
From (\ref{summary}) we therefore obtain a commutative diagram of perfect bilinear maps
\[
\xymatrix{
\bar T(G)/NT(G)\ar@{->>}[d]&*-<3pc>{\times}& C_G(N,\bar T(G))\ar[rr]^{\qquad\langle\cdot,\cdot\rangle_{N,\bar T(G)}^C}&&\dbZ/p\ar@{=}[d]\\
\bar T(G)/\pi\inv[T(G/N)]&*-<3pc>{\times}& B_G(N,\bar T(G))\ar@{^{(}->}[u]\ar[rr]^{\qquad\langle\cdot,\cdot\rangle_{N,\bar T(G)}^B}&&\dbZ/p.
}
\]
Hence the left vertical map is an equality if and only if the right vertical map is an equality.
This gives the required equivalence.
\end{proof}

\section{Unipotent upper-triangular matrices}
\label{section on unitriangular matrices}
We fix an integer $n\geq2$.
For a commutative unital profinite ring $R$ let $\dbU_n(R)$ be the group of all unipotent upper-triangular $(n+1)\times(n+1)$-matrices over $R$.
The additive group $R^+$ of $R$ embeds as a central subgroup of $\dbU_n(R)$ via  $r\mapsto I_{n+1}+rE_{1,n+1}$, where $I_{n+1}$ is the identity matrix and $E_{1,n+1}$ is the matrix with $1$ at entry $(1,n+1)$ and $0$ elsewhere.
Setting $\bar\dbU_n(R)=\dbU_n(R)/R^+$ we obtain a central extension
\[
\omega\colon \qquad 0\to R^+\to\dbU_n(R)\to\bar\dbU_n(R)\to 1.
\]

We define maps $\gam_1,\gam_2\colon\bar\dbU_n(R)\to\dbU_n(R)$ by mapping a coset $\bar M$ of $M\in\dbU_n(R)$ to the matrix in $\dbU_n(R)$ which is $0$ at entries $(1,j)$, $j=2,3\nek n+1$ (resp., $(i,n+1)$, $i=1,2\nek n$), and which coincides with $M$ elsewhere.
It is straightforward to verify that $\gam_1,\gam_2$ are homomorphisms, and $\Ker(\gam_1)\cap\Ker(\gam_2)=\{I_{n+1}\}$.
Therefore assumption (I) of \S\ref{section on T bar T} is satisfied for $\omega$.
Assumption (II) is satisfied when $p||R|<\infty$.

\begin{lem}
\label{bar T for Zassenhaus}
For every profinite group $G$ one has $T^{\bar\dbU_n(R)}(G)=T^{\dbU_{n-1}(R)}(G)$.
\end{lem}
\begin{proof}
There is an embedding $\dbU_{n-1}(R)\hookrightarrow\bar\dbU_n(R)$ into the upper-left $n\times n$-block, and with zeros at entries $(i,n+1)$, $i=2\nek n$.
By Lemma \ref{intersection of kernels}, it implies that $T^{\bar\dbU_n(R)}(G)\leq T^{\dbU_{n-1}(R)}(G)$.
Moreover, applying Lemma \ref{intersection of kernels} for the projections $\bar\dbU_n(R)\to\dbU_{n-1}(R)$ on the upper-left and lower-right $n\times n$ blocks shows that $T^{\dbU_{n-1}(R)}(G)\leq T^{\bar\dbU_n(R)}(G)$.
\end{proof}

Now suppose that the ring $R$ is finite.
Let $\bar G$ be a profinite group.
Given a continuous homomorphism $\bar\rho\colon \bar G\to\bar\dbU_n(R)$, we write $\bar\rho_{ij}$ for the $(i,j)$-entry of $\bar\rho$.
The maps $\bar\rho_{i,i+1}\colon\bar G\to R^+$ are group homomorphisms.
The pullbacks $\bar\rho^*(\alp_\omega)$ are the elements of the  \textsl{$n$-fold Massey products} $\langle\varphi_1\nek\varphi_n\rangle\subseteq H^2(\bar G,R^+)$, where $\varphi_1\nek\varphi_n\in H^1(\bar G,R^+)=\Hom_{\rm cont}(\bar G,R^+)$ \cite{Dwyer75}*{Th.\ 2.4}.
Specifically, $\langle\varphi_1\nek\varphi_n\rangle$ consists of all pullbacks $\bar\rho^*(\alp_\omega)$ such that $\varphi_i=\bar\rho_{i,i+1}$, $i=1,2\nek n$.
For our purposes, this can be taken as the definition of the Massey product in the profinite group setting.

In particular, when $n=2$ one has  $\langle\varphi_1,\varphi_2\rangle=\{\varphi_1\cup\varphi_2\}$, where the cup product is induced from the product map $R\tensor R\to R$.
Hence $\bar\rho^*(\alp_\omega)=\bar\rho_{12}\cup\bar\rho_{23}$ \cite{EfratMinac11}*{Prop.\ 9.1}.
Since any $\varphi_1,\varphi_2\in H^1(\bar G,R^+)$ can be realized as $\bar\rho_{12},\bar\rho_{23}$, respectively, for some $\bar\rho$, we obtain in this way all cup products $\varphi_1\cup\varphi_2$ in $H^2(\bar G,R^+)$.

\section{The $p$-Zassenhaus filtration}
\label{section on p Zassenhaus filtration}
We now focus on the ring $R=\dbZ/p$, with $p$ prime.
We abbreviate
\[
\dbU_n=\dbU_n(\dbZ/p), \quad \bar\dbU_n=\bar\dbU_n(\dbZ/p).
\]
Let $\Omega=\{\omega\}$, where $\omega=\omega_n$ is the central extension
\[
0\to\dbZ/p\to\dbU_n\to\bar\dbU_n\to1.
\]

For a profinite group $G$ we have, by Lemma \ref{bar T for Zassenhaus},
\[
T(G)=T^{\dbU_n}(G), \quad \bar T(G)=T^{\bar \dbU_n}(G)=T^{\dbU_{n-1}}(G).
\]

We set $\bar G=G/\bar T(G)=G/T^{\dbU_{n-1}}(G)$ and let $\pi\colon G\to \bar G$ be the natural projection.
In view of the remarks in \S\ref{section on unitriangular matrices}, we denote the subgroup of $H^2(\bar G)$ generated by all $\bar\rho^*(\alp_\omega)$, for continuous homomorphisms $\bar\rho\colon \bar G\to \bar\dbU_n$, by $H^2(\bar G)_{n-{\rm Massey}}$.
Recall that $H^2(\bar G)_\pi$ is the subgroup of $H^2(\bar G)$ generated by all such pullbacks which become trivial under the inflation to $H^2(G)$.
We denote it in the current setting by $H^2(\bar G)_{n-{\rm Massey},G}$.

Now Proposition \ref{a special case} gives:

\begin{thm}
\label{exact sequence and pairing for Zassenhaus filtration}
\begin{enumerate}
\item[(a)]
The bilinear map $\langle\cdot,\cdot\rangle^B$  gives a perfect bilinear map
\[
T^{\dbU_{n-1}}(G)/T^{\dbU_n}(G)\times H^2(G/T^{\dbU_{n-1}}(G))_{n-{\rm Massey},G}\to\dbZ/p,
\]
\item[(b)]
There is an exact sequence
\[
0\to H^2(G/T^{\dbU_{n-1}}(G))_{n-{\rm Massey},G}\hookrightarrow H^2(G/T^{\dbU_{n-1}}(G))\xrightarrow{\inf}H^2(G/T^{\dbU_n}(G)).
\]
\end{enumerate}
\end{thm}

Corollary \ref{generalization of a result of Dwyer} gives in this setting:

\begin{prop}
\label{isomorphic kernels for the Zassenhaus filtration}
Let $N_1\leq N_2$ be closed normal subgroups of $G$ contained in $T^{\dbU_{n-1}}(G)$.
The following conditions are equivalent:
\begin{enumerate}
\item[(a)]
$H^2(G/N_i)_{n-{\rm Massey},G}$, $i=1,2$, are isomorphic via inflation;
\item[(b)]
$N_1T^{\dbU_n}(G)= N_2T^{\dbU_n}(G)$.
\end{enumerate}
\end{prop}

The above structural results can sometimes be rephrased in terms of the \textsl{$p$-Zassenhaus filtration} $G_{(n,p)}$, $n=1,2,3\nek$ of $G$  (see the Introduction).
By its definition, the quotients $G_{(n-1,p)}/G_{(n,p)}$ are elementary abelian $p$-groups, so by induction, $G/G_{(n,p)}$ is a pro-$p$ group for every $n$.

When $S$ is a free profinite group it was shown in \cite{Efrat14a}*{Th.\ A} that
\begin{equation}
\label{TS as Zassenhaus subgroup}
T(S)=T^{\dbU_n}(S)=S_{(n+1,p)}.
\end{equation}
Consequently, $ \bar T(S)=T^{\dbU_{n-1}}(S)=S_{(n,p)}$.
See also \cite{MinacTan15}*{Th.\ 2.7} and \cite{Efrat14b} for alternative proofs.

Therefore Theorem \ref{exact sequence and pairing for Zassenhaus filtration} and  Proposition \ref{isomorphic kernels for the Zassenhaus filtration}, when applied to $G=S$, recover \cite{Efrat14a}*{Th.\ B}
and \cite{Efrat14a}*{Cor.\ 10.2}, respectively.

Now Theorem \ref{The transfer theorem general form} can be interpreted in this setting as follows:

\begin{thm}
\label{equivalence for the Zassenhaus filtration}
Let $S$ be a free profinite group,  let $N$ be a closed normal subgroup of $S$ contained in $S_{(n,p)}$, and let $G=S/N$.
Then:
\begin{enumerate}
\item[(a)]
$T^{\dbU_{n-1}}(G)=G_{(n,p)}$.
\item[(b)]
One has $T^{\dbU_n}(G)=G_{(n+1,p)}$ if and only if
\[
\Ker\bigl(H^2(G/G_{(n,p)})_{n-{\rm Massey}}\xrightarrow{\inf}H^2(G)\bigr)
\]
is generated by elements of $n$-fold Massey products.
\end{enumerate}
\end{thm}
\begin{proof}
Let $\pi\colon S\to G$ be the natural projection.

\medskip

(a) \quad
As $N\leq \bar  T(S)$, Lemma \ref{functoriality of T, bar T} gives:
\[
T^{\dbU_{n-1}}(G)=\bar T(G)=\pi[\bar T(S)]=\pi[T^{\dbU_{n-1}}(S)]=\pi[S_{(n,p)}]=G_{(n,p)}.
\]

(b) \quad
One has $S/S_{(n,p)}\isom G/G_{(n,p)}$ via $\pi$.
Hence the condition about the kernel in (b) means that $\Ker(H^2(S/S_{(n,p)})_{n-{\rm Massey}}\to H^2(S/N))$ is generated by $n$-fold Massey products.
This is exactly the kernel generating condition for the subgroups $N,S_{(n,p)}$ of $S$ (see Remark \ref{remarks on liftable}(3)).
By Theorem \ref{The transfer theorem general form} (with $G$ replaced by $S$), it is equivalent to $T(G)=\pi[T(S)]$, i.e., to $T^{\dbU_n}(G)=G_{(n+1,p)}$.
\end{proof}

The condition about the kernel in Theorem  \ref{equivalence for the Zassenhaus filtration}(b) is called in \cite{Efrat14a} the \textsl{$n$-Massey kernel condition} on $G$.
Thus the ``if" part in the equivalence recovers \cite{Efrat14a}*{Th.\ A'}.

\begin{rem}
\label{G and G(p)}
\rm
Denote the maximal pro-$p$ quotient of $G$ by $G(p)$.
Let $\lam\colon G\to G(p)$ be the natural epimorphism.
Then $\Ker(\lam)$ is contained in both $T(G)=T^{\dbU_n(\dbZ/p)}(G)$ and $G_{(n+1,p)}$.
By Lemma \ref{functoriality of T, bar T}, $\lam[T(G)]=T(G(p))$.
Also, $\lam[G_{(n+1,p)}]=G(p)_{(n+1,p)}$.
Consequently, $T(G)=G_{(n+1,p)}$ if and only $T(G(p))=(G(p))_{(n+1,p)}$ (and similarly for $\bar T(G)$, $G_{(n,p)}$).
\end{rem}

\section{The lower $p$-central filtration}
\label{section on lower p-central extension}
We fix $n\geq2$.
For every $1\leq s\leq n$ let $\dbU_{n,s}=\dbU_s(\dbZ/p^{n-s+1})$ and $\bar\dbU_{n,s}=\bar\dbU_{n,s}(\dbZ/p^{n-s+1})$ (with notation as in \S\ref{section on unitriangular matrices}).
Let $\Omega$ consist of the central extensions
\[
\omega_{n,s}\colon\quad
0\to\dbZ/p\to\dbU_{n,s}\to\bar\dbU_{n,s}\to1,\quad s=1,2\nek n.
\]
For example, the extensions $\omega_{n,1},\omega_{n,n}$ are
\[
0\to\dbZ/p\to\dbZ/p^n\to\dbZ/p^{n-1}\to0, \quad  0\to \dbZ/p\to\dbU_n(\dbZ/p)\to\bar\dbU_n(\dbZ/p)\to 1,
\]
respectively.
As noted in \S\ref{section on unitriangular matrices}, assumptions (I) and (II) are satisfied.

For every profinite group $G$ we clearly have
\[
 T(G)=\bigcap _{s=1}^nT^{\dbU_{n,s}}(G).
\]

\begin{prop}
\label{T bar in ut case}
One has
$\bar T(G)=\bigcap _{s=1}^{n-1}T^{\dbU_{n-1,s}}(G)$.
\end{prop}
\begin{proof}
For every $2\leq s\leq n$, the group $\dbU_{n-1,s-1}$ embeds as the upper-left $s\times s$-block of $\bar\dbU_{n,s}$, with zeros at entries $(i,s+1)$, $i=2,3\nek s$.
By Lemma \ref{intersection of kernels},
\[
T^{\bar \dbU_{n,s}}(G)\leq T^{\dbU_{n-1,s-1}}(G), \quad s=2\nek n.
\]
Moreover, $\bar\dbU_{n,1}\isom\dbZ/p^{n-1}\isom\dbU_{n-1,1}$, so
\[
T^{\bar\dbU_{n,1}}(G)=T^{\dbU_{n-1,1}}(G).
\]
Hence
\[
\bar T(G)=\bigcap_{s=1}^nT^{\bar\dbU_{n,s}}(G)\leq\bigcap _{s=1}^{n-1}T^{\dbU_{n-1,s}}(G).
\]

On the other hand, for $2\leq s\leq n-1$ let $\gam_1,\gam_2\colon\bar\dbU_{n,s}\to\dbU_{n-1,s-1}$ be the projections on the upper-left and lower-right $s\times s$-blocks.
Also let $\gam_3\colon\bar\dbU_{n,s}\to\dbU_{n-1,s}$ be the homomorphism induced by the ring epimorphism $\dbZ/p^{n-s+1}\to\dbZ/p^{n-s}$.
Then $\Ker(\gam_1)\cap\Ker(\gam_2)\cap\Ker(\gam_3)=\{I_s\}$.
Lemma \ref{intersection of kernels} therefore implies that
\[
T^{\dbU_{n-1,s-1}}(G)\cap T^{\dbU_{n-1,s}}(G)\leq T^{\bar\dbU_{n,s}}(G), \quad s=2\nek n-1.
\]
Moreover, Lemma \ref{bar T for Zassenhaus} gives
$T^{\dbU_{n-1,n-1}}(G)=T^{\bar \dbU_{n,n}}(G)$.
We deduce that
\[
\bigcap_{s=1}^{n-1}T^{\dbU_{n-1,s}}(G)\leq\bigcap_{s=1}^nT^{\bar\dbU_{n,s}}(G)=\bar T(G),
\]
and the assertion follows.
\end{proof}

When $S$ is a free profinite group one has
\[
T(S)=S^{(n+1,p)}.
\]
Indeed, this is proved in \cite{MinacTan15}*{Th.\ 2.7} when $S$ is finitely generated; The general case follows by an inverse limit argument.
In view of Proposition \ref{T bar in ut case}, this implies that $\bar T(S)=S^{(n,p)}$.
See \cite{Efrat14b} for a variant of this result.

Given a profinite group $G$, we call the pullbacks $\bar\rho^*(\alp_{\omega_{n,s}})$, where $1\leq s\leq n$ and $\bar\rho\colon G/\bar T(G)\to\bar\dbU_{n,s}$ is a continuous homomorphism, \textsl{$n$-unitriangular elements}.
They were studied in \cite{Efrat17} under the name \textsl{the unitriangular spectrum}.
Let $H^2(G/\bar T(G))_{n-{\rm ut}}$ be the subgroup of $H^2(G/\bar T(G))$ generated by all $n$-unitriangular, and let $H^2(G/\bar T(G))_{n-{\rm ut},G}$ be the subgroup of $H^2(G/\bar T(G))$ generated by all these pullbacks with $\bar\rho$ liftable to $G$.

Let $G^{(n,p)}$, $n=1,2,3\nek$ be the \textsl{lower $p$-central filtration} of the profinite group $G$ (see the Introduction).
The quotients $G^{(n,p)}/G^{(n+1,p)}$ are elementary $p$-groups, so by induction, $G/G^{(n,p)}$ is a pro-$p$ group for every $n$.

\begin{prop}
\label{H2 n-ut is all of H2}
For a free profinite group $S$ one has
\[
H^2(S/\bar T(S))_{n-{\rm ut},S}=H^2(S/\bar T(S))_{n-{\rm ut}}=H^2(S/S^{(n,p)}).
\]
\end{prop}
\begin{proof}
The left equality follows from Remark \ref{remarks on liftable}(3).

The main result of \cite{Efrat17} gives a canonical linear basis (called the \textsl{Lyndon basis}) for $H^2(S/S^{(n,p)})$ which is made of certain $n$-unitiangular elements $\bar\rho^*(\alp_{\omega_{n,s}})$.
This gives the right equality.
\end{proof}

We deduce the analog of Theorem \ref{equivalence for the Zassenhaus filtration}:

\begin{thm}
\label{equivalence for lower p central filtration}
Let $S$ be a free profinite group,  let $N$ be a closed normal subgroup of $S$ contained in  $S^{(n,p)}$, and let $G=S/N$.
Then:
\begin{enumerate}
\item[(a)]
$\bar T(G)=G^{(n,p)}$.
\item[(b)]
One has $T(G)=G^{(n+1,p)}$ if and only if
\[
\Ker\bigl(H^2(G/G^{(n,p)})\xrightarrow{\inf}H^2(G)\bigr)
\]
is generated by $n$-unitriangular elements.
\end{enumerate}
\end{thm}
\begin{proof}
Let $\pi\colon S\to G$ be the natural projection.

\medskip

(a) \quad
As $N\leq \bar  T(S)$, Lemma \ref{functoriality of T, bar T} gives:
\[
\bar T(G)=\pi[\bar T(S)]=\pi[S^{(n,p)}]=G^{(n,p)}.
\]

(b) \quad
We have $S/S^{(n,p)}\isom G/G^{(n,p)}$ via $\pi$.
Hence the condition about the kernel means that $\Ker(\inf\colon H^2(S/S^{(n,p)})\to H^2(S/N))$ is generated by $n$-triangular elements.
In view of Proposition \ref{H2 n-ut is all of H2}, this is exactly the kernel generating condition for the subgroups $N,S^{(n,p)}$ of $S$, and $\Omega$ as above.
By Theorem \ref{The transfer theorem general form} (with $G$ replaced by $S$), it is equivalent to  $T(G)=\pi[T(S)]$, i.e., to $T(G)=\pi[S^{(n+1,p)}]=G^{(n+1,p)}$.
\end{proof}

When the condition on the kernel in (b) holds, we say that $G$ satisfies the \textsl{$n$-unitriangular kernel condition}.

From Theorem \ref{equivalence for lower p central filtration} and Proposition  \ref{a special case} we deduce:

\begin{cor}
Let $S$ be a free profinite group,  let $N$ be a closed normal subgroup of $S$ contained in  $S^{(n,p)}$, and suppose that $G=S/N$ satisfies the $n$-unitriangular kernel condition.
Let  $\pi\colon G\to G/G^{(n,p)}$ be the natural projection.
Then:
\begin{enumerate}
\item[(a)]
There is a non-degenerate bilinear map
\[
G^{(n,p)}/G^{(n+1,p)}\times H^2(G/G^{(n,p)})_\pi\to\dbZ/p.
\]
\item[(b)]
There is an exact sequence
\[
0\to H^2(G/G^{(n,p)})_\pi\to H^2(G/G^{(n,p)})\xrightarrow{\inf}H^2(G/G^{(n+1,p)}).
\]
\end{enumerate}
\end{cor}

\begin{rem}
\rm
Similarly to Remark \ref{G and G(p)}, we have $T(G)=G^{(n+1,p)}$ if and only $T(G(p))=(G(p))^{(n+1,p)}$, and likewise for $\bar T(G)$, $G^{(n,p)}$.
\end{rem}

\section{The lower $p$-central filtration for $p>2$, $n=2$}
\label{section on the lower central filtration for odd prime}
The discussion in this section is inspired by  \cite{EfratMinac11}.
Let $p>2$ and let $M_{p^3}$ be the extra-special group of order $p^3$ and exponent $p^2$.
Thus
\begin{equation}
\label{presentation of Mp3}
M_{p^3}=\bigl\langle r,s\ \bigm|  r^{p^2}=s^p=1, [r,s]=r^p\bigr\rangle.
\end{equation}
Let $\Omega$ consist of the central extensions
\[
\omega\colon \  0\to\dbZ/p\to\dbZ/p^2\to \dbZ/p\to0, \quad \omega'\colon \ 0\to\dbZ/p\to M_{p^3}\to(\dbZ/p)^2\to0.
\]
where the epimorphism in $\omega'$ is given by $r\mapsto(1,0)$, $s\mapsto(0,1)$.
The group $\dbZ/p$ clearly embeds in $\dbZ/p^2$, and the group $(\dbZ/p)^2$ embeds in $M_{p^3}$ via $(1,0)\mapsto r^p$, $(0,1)\mapsto sr^p$ \cite{EfratMinac11}*{Remark 8.1(c)}.
Hence assumptions (I), (II) of \S\ref{section on T bar T} are satisfied.

For a profinite group $\bar G$, we consider $\omega$ as a short exact sequence of trivial $\bar G$-modules.
Then the \textsl{Bockstein homomorphism}  $\Bock=\Bock_{\bar G}\colon H^1(\bar G)\to H^2(\bar G)$ is the connecting homomorphism in the associated long exact sequence of cohomology groups \cite{NeukirchSchmidtWingberg}*{Th.\ 1.3.2}.

\begin{prop}
\label{pullbacks in mixed case}
\begin{enumerate}
\item[(a)]
For a continuous homomorphism $\bar\rho\colon \bar G\to\dbZ/p$ we have $\bar\rho^*(\alp_{\omega})=\Bock(\bar\rho)$.
\item[(b)]
For a continuous epimorphism $\bar\rho=(\bar\rho_1,\bar\rho_2)\colon\bar G\to(\dbZ/p)^2$ we have $\bar\rho^*(\alp_{\omega'})=\Bock(\bar\rho_1)+\bar\rho_1\cup\bar\rho_2$.
\item[(c)]
For a continuous homomorphism $\bar\rho\colon\bar G\to(\dbZ/p)^2$ which is not surjective we have either  $\bar\rho^*(\alp_{\omega'})=\Bock(\bar\rho)$ (where $\bar\rho$ is considered as an element of $H^1(\bar G)$), or  $\bar\rho^*(\alp_{\omega'})=0$.
\end{enumerate}
\end{prop}
\begin{proof}
(a) \quad
See \cite{EfratMinac11}*{Prop.\ 9.2} for the case $\bar\rho\neq0$.
The case $\bar\rho=0$ is trivial.
\medskip

(b) \quad
See \cite{EfratMinac11}*{Prop.\ 9.4}.

\medskip

(c) \quad
We may assume that $\bar\rho\neq0$, so $\Img(\bar\rho)\isom\dbZ/p$.
Then $\bar\rho$ breaks via a commutative diagram
\[
\xymatrix{
&&&\bar G\ar@{->>}[d]^{\bar\rho}&\\
0\ar[r]&\dbZ/p\ar@{=}[d]\ar[r]&D\ar[r]\ar@{_{(}->}[d]&\dbZ/p\ar[r]\ar@{_{(}->}[d]&0\\
0\ar[r]&\dbZ/p\ar[r]& M_{p^3}\ar[r]&(\dbZ/p)^2\ar[r]&0,
}
\]
where $D$ is either $\dbZ/p^2$ or $(\dbZ/p)^2$.
In the first case,  $\bar\rho^*(\alp_{\omega'})=\Bock(\bar\rho)$.
In the second case, the upper extension splits, so $\bar\rho^*(\alp_{\omega'})=0$.
\end{proof}

Now let $G$ be a profinite group and set $\bar G=G/G^{(2,p)}$.

\begin{lem}
\label{dichotomy}
Let $\psi,\xi\colon G\to\dbZ/p$ be continuous homomorphisms such that $\Bock_G(\psi)=\psi\cup\xi$.
Let $\bar\psi,\bar\xi\colon \bar G\to\dbZ/p$ be the homomorphisms induced by $\psi,\xi$, respectively.
Then
\begin{enumerate}
\item[(i)]
$\bar\psi^*(\alp_{\omega})$ is $G$-liftable; or
\item[(ii)]
$\bar\rho^*(\alp_{\omega'})$ is $G$-liftable, where $\bar\rho=(\bar\psi,-\bar\xi)$.
\end{enumerate}
\end{lem}
\begin{proof}
(Compare \cite{EfratMinac11}*{Prop.\ 10.2}).
First assume that $\Bock_G(\psi)=0$.
Since $\bar\psi^*(\alp_{\omega})=\Bock_{\bar G}(\bar\psi)$ (Proposition \ref{pullbacks in mixed case}(a)), the functoriality of $\Bock$ (as a connecting homomorphism) implies that $\inf_G(\bar\psi^*(\alp_{\omega}))=\Bock_G(\psi)=0$, i.e., $\bar\psi^*(\alp_{\omega})$ is $G$-liftable, and (i) holds.

Next suppose that $\Bock_G(\psi)\neq0$.
Then $\psi\cup\xi\neq0$.
As $p>2$, $\psi\cup\psi=0$, so $\psi,\xi$ are $\dbF_p$-linearly independent in $H^1(G)$.
Hence $\bar\psi,\bar\xi$ are $\dbF_p$-linearly independent in $H^1(\bar G)$.
Therefore $\bar\rho=(\bar\psi,-\bar\xi)\colon \bar G\to(\dbZ/p)^2$ is an epimorphism.
One has $\bar\rho^*(\alp_{\omega'})=\Bock_{\bar G}(\bar\psi)-\bar\psi\cup\bar\xi$ (Proposition \ref{pullbacks in mixed case}(b)), so by hypothesis, $\inf_G(\bar\rho^*(\alp_{\omega'}))=0$, whence (ii).
\end{proof}

Let $\pi\colon G\to \bar G$ be the natural projection.
By Proposition \ref{pullbacks in mixed case}, $H^2(\bar G)_\pi$ is the subgroup of $H^2(\bar G)$ generated by the elements of the forms $\Bock(\bar\rho)$ and $\Bock(\bar\rho)+\bar\rho\cup\bar\rho'$ which vanish under inflation to $H^2(G)$.

The quotient $\bar G$ is an elementary abelian $p$-group, so $H^2(\bar G)$ is generated by all Bockstein elements $\Bock(\bar\psi)$ and all cup products $\bar\psi\cup\bar\psi'$, with $\bar\psi,\bar\psi'\in H^1(\bar G)$ \cite{EfratMinac11}*{Cor.\ 2.9(a)}.
Therefore it is generated by the cohomology elements $\Bock(\bar\rho)$ and $\Bock(\bar\rho)+\bar\rho\cup\bar\rho'$.

We say that $G$ satisfies the \textsl{Bockstein-cup kernel condition} if
\[
\Ker\bigl(H^2(\bar G)\xrightarrow{\inf}H^2(G)\bigr)
\]
is generated by elements of the forms $\Bock(\bar\rho)$ and $\Bock(\bar\rho)+\bar\rho\cup\bar\rho'$.

Next we note that
\[
\bar T(G)=T^{\dbZ/p}(G)=G^{(2,p)}, \qquad T(G)=T^{\dbZ/p^2}(G)\cap T^{M_{p^3}}(G).
\]
Since $(\dbZ/p^2)^{(3,p)}=1$ and $(M_{p^3})^{(3,p)}=1$  \cite{EfratMinac11}*{Remark 8.1(b)}, we deduce that
\[
G^{(3,p)}\leq T(G)\leq G^{(2,p)}.
\]

\begin{lem}
\label{S3p=TS}
If $S$ is a free profinite group, then $S^{(3,p)}=T(S)$.
\end{lem}
\begin{proof}
By Proposition \ref{transgression perfect pairing},  $\langle\cdot,\cdot\rangle^A_{1,S^{2,p)}}$ gives a perfect bilinear map as in the upper row of the following commutative diagram:
\[
\xymatrix{
S^{(2,p)}/S^{(3,p)}\ar@{->>}[d]&*-<3pc>{\times}& H^2(S/S^{(2,p)})\ar[r]&\dbZ/p\ar@{=}[d]\\
S^{(2,p)}/T(S)&*-<3pc>{\times}&H^2(S/S^{(2,p)})_\pi\ar@{^{(}->}[u]\ar[r]&\dbZ/p.
}
\]
The lower row of the diagram is perfect by Proposition \ref{a special case}(a).
By Proposition \ref{H2 n-ut is all of H2}, the right vertical map is an equality.
Therefore the left vertical map is also an equality.
\end{proof}

We extend Lemma \ref{S3p=TS} as follows:

\begin{thm}
\label{equivalence for Bockstein cup context}
Let $S$ be a free profinite group,  let $N$ be a closed normal subgroup of $S$ contained in $S^{(2,p)}$, and let $G=S/N$.
Then $T(G)=G^{(3,p)}$ if and only if $G$ satisfies the Bockstein-cup kernel condition.
\end{thm}
\begin{proof}
We have $S/S^{(2,p)}\isom G/G^{(2,p)}=\bar G$ via $\pi$.
Hence the Bockstein-cup kernel condition means that $\Ker(\inf\colon H^2(S/S^{(2,p)})\to H^2(S/N))$ is generated by elements of the forms $\Bock(\bar\rho)$ and $\Bock(\bar\rho)+\bar\rho\cup\bar\rho'$.
In view of Remark \ref{remarks on liftable}(3), this is exactly the kernel generating condition for the subgroups $N,S^{(2,p)}$ of $S$.
By Theorem \ref{The transfer theorem general form} (with $G$ replaced by $S$), it is equivalent to $T(G)=\pi[T(S)]$.
In light of Lemma \ref{S3p=TS}, the latter equality means that $T(G)=\pi[S^{(3,p)}]=G^{(3,p)}$.
\end{proof}

Let $H^2(G)_{\rm cup,Bock}$ be the subgroup of $H^2(G)$ generated by cup products and Bockstein elements.
Then Proposition \ref{a special case} gives:

\begin{cor}
Let $G$ be a profinite group satisfying  the Bockstein-cup kernel condition.
Then:
\begin{enumerate}
\item[(a)]
There is a canonical perfect bilinear map
\[
G^{(2,p)}/G^{(3,p)}\times H^2(\bar G)_{\rm cup,Bock}\to\dbZ/p,
\]
\item[(b)]
There is an exact sequence
\[
0\to H^2(\bar G)_{\rm cup,Bock}\to H^2(\bar G)\to H^2(G/G^{(3,p)}).
\]
\end{enumerate}
\end{cor}
\begin{proof}
Take a free profinite group $S$ and a continuous epimorphism $\pi\colon S\to G$ whose kernel $N$ is contained in the Frattini subgroup of $S$.
In particular, $N\leq S^{(2,p)}$, so by Theorem \ref{equivalence for Bockstein cup context}, $T(G)=G^{(3,p)}$.
We now apply Proposition \ref{a special case} with $G$ replaced by $S$.
\end{proof}

\section{Applications to absolute Galois groups}
\label{section on absolute Galois groups}
The next theorem assembles the facts in Galois cohomology which are needed for the intersection theorems for absolute Galois groups of \cite{MinacSpira96}, \cite{EfratMinac11}, and \cite{EfratMinac 17}, given in the Introduction.
It is based on the injectivity part of the Merkurjev--Suslin theorem (\cite{MerkurjevSuslin82}, \cite{GilleSzamuely06}),  the identification of Bockstein elements as cup products \cite{EfratMinac11}*{Prop.\ 2.6}, as well as the structure of $H^2$ for elementary abelian $p$-groups.
While all ingredients of the proof essentially appear in the above-mentioned references, our new formalism demonstrates the common methodology which underlies all the previous proofs.

For a field $F$ containing a root of unity $\zeta$ of order $p$, whence the full group $\mu_p$ of $p$th roots of unity, we identify $\mu_p=\dbZ/p$ as $G_F$-modules, where $\zeta$ corresponds to $\bar1\in\dbZ/p$.
This gives the Kummer isomorphism $F^\times/(F^\times)^p\xrightarrow{\sim}H^1(G_F)$.
Given $a\in F^\times$, let $(a)_F\in H^1(G_F)$ correspond to the coset of $a$ in $F^\times/(F^\times)^p$.

By \cite{EfratMinac11}*{Prop.\ 2.6}, for every $\varphi\in H^1(G_F)$ one has
\begin{equation}
\label{Bockstein identity}
\Bock_{G_F}(\varphi)=\varphi\cup(\zeta)_F.
\end{equation}

\begin{thm}
\label{MS}
Let $F$ be a field containing a root of unity $\zeta$ of order $p$, and let $G=G_F$ be its absolute Galois group.
Then:
\begin{enumerate}
\item[(a)]
$G$ satisfies the $2$-Massey kernel condition;
\item [(b)]
When $p=2$,  $G$ satisfies the $2$-unitriangular kernel condition;
\item[(c)]
When $p>2$, $G$ satisfies the Bockstein--cup kernel condition.
\end{enumerate}
\end{thm}
\begin{proof}
We note that $G_{(2,p)}=G^{(2,p)}$.
Set $\bar G=G/G_{(2,p)}=G/G^{(2,p)}$.
Then $H^1(\bar G)\isom H^1(G)$ via inflation.
We recall that the 2-fold Massey product is just the usual cup product map $\cup\colon H^1(\bar G)^{\tensor2}\to H^2(\bar G)$
(Its image $H^2(\bar G)_{2-{\rm Massey}}$ it is also denoted by $H^2(\bar G)_{\rm dec}$, for the \textsl{decomposable} part of $H^2(G)$ \cite{CheboluEfratMinac12}).

\medskip

(a) \quad
By the injectivity part of the Merkurjev--Suslin theorem, the kernel of the cup product map $\cup\colon H^1(G)^{\tensor2}\to H^2(G)$ is the \textsl{``Steinberg group"},
generated by all tensor products $(a)_F\tensor(1-a)_F$, with $a\in F$, $a\neq0,1$.
Therefore this kernel is also generated by all $\varphi_1\tensor\varphi_2$, where $\varphi_1,\varphi_2\in H^1(G) $ satisfy $\varphi_1\cup\varphi_2=0$.
There is a commutative square
\[
\xymatrix{
H^1(\bar G)^{\tensor2}\ar[r]^{\inf}_{\sim}\ar@{->>}[d]_{\cup}& H^1(G)^{\tensor2}\ar[d]^{\cup}\\
H^2(\bar G)_{2-{\rm Massey}}\ar[r]^{\quad\inf}&H^2(G).
}
\]
It follows that the kernel of the lower inflation map is generated by cup products from $H^1(\bar G)^{\tensor2}$ which are in this kernel, as desired.

\medskip

(b)\quad
We recall that for every prime $p$, $H^2(\bar G)_{2-{\rm ut}}$ is the subgroup of $H^2(\bar G)$ generated by all Bockstein elements and all cup products (see Proposition \ref{pullbacks in mixed case}(a)).
As $p=2$,  $\bar G$ is an elementary $2$-abelian group, so $H^2(\bar G)$ is generated by cup products only \cite{EfratMinac11}*{Cor.\ 2.9(b)}.
Therefore $H^2(\bar G)_{2-{\rm ut}}=H^2(\bar G)_{2-{\rm Massey}}=H^2(\bar G)$.
The assertion therefore follows from (a).

\medskip

(c) \quad
Consider $\alp$ in $\Ker\bigl(\inf\colon H^2(\bar G)\to H^2(G)\bigr)$.
It is a sum of Bockstein elements and cup products in $H^2(\bar G)$.
Therefore it can be written as
\begin{equation}
\label{alpha as a sum}
\alp=\sum_i\bar\psi_i\cup\bar\psi'_i+(\Bock_{\bar G}(\bar\varphi)-\bar\varphi\cup\overline{(\zeta)}_F),
\end{equation}
with $\bar\psi_i,\bar\psi'_i,\bar\varphi\in H^1(\bar G)$, and where $\overline{(\zeta)}_F\in H^1(\bar G)$ corresponds to the Kummer element $(\zeta)_F\in H^1(G)$ under the inflation isomorphism $H^1(\bar G)\isom H^1(G)$.
By (\ref{Bockstein identity}), $\Bock_{\bar G}(\bar\varphi)-\bar\varphi\cup\overline{(\zeta)}_F$ is always in the above kernel.
Hence $\sum_i\bar\psi_i\cup\bar\psi'_i$ is also in the kernel.
By (a), we may replace it by another such sum, to assume without loss of generality that each $\bar\psi_i\cup\bar\psi'_i$ separately is in this kernel.
Moreover, $\Bock_{\bar G}(\bar\varphi)-\bar\varphi\cup\overline{(\zeta)}_F=\Bock_{\bar G}(\bar\varphi)+\bar\varphi\cup(-\overline{(\zeta)}_F)$, is a pullback corresponding to the extension $\omega'$ of \S\ref{section on the lower central filtration for odd prime} (Proposition \ref{pullbacks in mixed case}(b)).
Then in (\ref{alpha as a sum}) all the summands are pullbacks in the kernel, as desired.
\end{proof}

We now recover in a uniform way the intersection theorems by Min\'a\v c, Spira, and the author (\cite{EfratMinac 17}*{Th.\ D},  \cite{MinacSpira96}*{Cor.\ 2.18}, \cite{EfratMinac11}).

\begin{thm}
\label{Efrat Minac Spira}
Let $G=G_F$ be the absolute Galois group of a field $F$ which contains a root of unity  $\zeta$ of order $p$.
Then:
\begin{enumerate}
\item[(a)]
$G_{(3,p)}$ is the intersection of all closed normal subgroups $M$ of $G$ such that $G/M$ is either $\{1\}$, $\dbZ/p$, or $\dbU_2(\dbZ/p)$.
\item[(b)]
When $p=2$, $G^{(3,2)}$ is the intersection of all closed normal subgroups $M$ of $G$ such that $G/M$ is isomorphic to either $\{1\}$, $\dbZ/2$, $\dbZ/4$, or $\dbU_2(\dbZ/2)=D_4$.
\item[(c)]
When $p>2$, $G^{(3,p)}$ is the intersection of all closed normal subgroups $M$ of $G$ such that $G/M$ is isomorphic to either $\{1\}$, $\dbZ/p^2$, or $M_{p^3}$.
\end{enumerate}
\end{thm}

Note that the quotient $\{1\}$ is needed only when $F$ is separably closed.

\begin{proof}
Take a presentation $G=S/N$ of $G$, with $S$ a free profinite group and $N$ contained in its Frattini subgroup.
Then $N\leq S_{(2,p)}=S^{(2,p)}$.
Let $\bar G=G/G_{(2,p)}=G/G^{(2,p)}$

\medskip

(a) \quad
Let $\Omega$ consist of the single extension
\[
0\to\dbZ/p\to \dbU_2(\dbZ/p)\to\bar\dbU_2(\dbZ/p)\to1.
\]
By Theorem \ref{MS}(a) and Theorem \ref{equivalence for the Zassenhaus filtration}(b), $G_{(3,p)}=T^{\dbU_2(\dbZ/p)}(G)$.
The subgroups of $\dbU_2(\dbZ/p)$ are $\{1\},\dbZ/p,(\dbZ/p)^2,\dbU_2(\dbZ/p)$.
Moreover, any closed normal subgroup of $G$ with quotient $(\dbZ/p)^2$ can be replaced in the intersection by two closed normal subgroups with quotient $\dbZ/p$.

\medskip

(b) \quad
Let $\Omega$ consist of the central extensions
\[
0\to\dbZ/2\to\dbZ/4\to\dbZ/2\to0, \quad 0\to\dbZ/2\to\dbU_2(\dbZ/2)=D_4\to(\dbZ/2)^2\to1.
\]
By Theorem \ref{MS}(b) and Theorem \ref{equivalence for lower p central filtration}(b),  $G^{(3,2)}=T(G)=T^{\dbZ/4}(G)\cap T^{D_4}(G)$.
The subgroups of $\dbZ/4$ and $D_4$ are $\{1\},\dbZ/2,\dbZ/4,(\dbZ/2)^2$, and $D_4$.
Again, a closed normal subgroup of $G$ with quotient $(\dbZ/2)^2$ can be replaced in the intersection by two closed normal subgroups with quotient $\dbZ/2$.

\medskip

(c) \quad
Let $\Omega$ consist of the central sequences
\[
\omega\colon \  0\to\dbZ/p\to\dbZ/p^2\to \dbZ/p\to0, \quad \omega'\colon \  0\to\dbZ/p\to M_{p^3}\to(\dbZ/p)^2\to0.
\]
By Theorem \ref{MS}(c) and Theorem \ref{equivalence for Bockstein cup context}, $G^{(3,p)}=T(G)=T^{\dbZ/p^2}(G)\cap T^{M_{p^3}}(G)$.
The subgroups of $\dbZ/p^2$ and $M_{p^3}$ are $\{1\}$, $\dbZ/p$, $\dbZ/p^2$, $(\dbZ/p)^2$ and $M_{p^3}$.
As before, the quotient $(\dbZ/p)^2$ can be omitted from the list.

It remains to show that the quotient $\dbZ/p$ can also be omitted from the list.
To this end, let $\psi\colon G\to\dbZ/p$ be a continuous epimorphism, and let $\bar\psi\colon \bar G\to\dbZ/p$ be the induced epimorphism.
By Lemma \ref{dichotomy} and (\ref{Bockstein identity}), one of the following holds:

\noindent
Case (i): \quad $\bar\psi^*(\alp_{\omega})$ is $G$-liftable.  \quad
Thus $\psi$ lifts to a continuous homomorphism $\rho\colon G\to\dbZ/p^2$.
The epimorphism in $\omega$ is a Frattini cover \cite{FriedJarden08}*{Def.\ 22.5.1}, so $\rho$ is surjective.
Consequently, $\Ker(\psi)$ contains a closed normal subgroup $N$ of $G$ with $G/N\isom\dbZ/p^2$.

\medskip

\noindent
Case (ii): \quad
$(\bar\psi,-\overline{(\zeta)}_F)^*(\alp_{\omega'})$ is $G$-liftable (where $\overline{(\zeta)}_F\in H^1(\bar G)$ is as in the previous proof).  \quad
Thus $(\psi,-(\zeta)_F)$ lifts to a continuous homomorphism $\rho\colon G\to M_{p^3}$.
The epimorphism $\omega'$ is also a Frattini cover, so $\rho$ is surjective.
Then $\psi$ is the composition of $\rho$ with the epimorphism $M_{p^3}\to\dbZ/p$, given by $r\mapsto 1$ and $s\mapsto 0$, so $\Ker(\psi)$ contains a closed normal subgroup $N$ of $G$ with $G/N\isom M_{p^3}$.

Consequently, in both cases, the normal subgroup $\Ker(\psi)$ of $G$ with quotient $\dbZ/p$ may be replaced by a normal closed subgroup $N$ with $G/N\isom\dbZ/p^2,M_{p^3}$.
\end{proof}

The \textsl{``Kernel $n$-Unipotent Conjecture"} of Min\'a\v c and T\^an \cite{MinacTan15} predicts that, when $G=G_F$ is the absolute Galois group of a field $F$ containing a root of unity of order $p$,  the equality $G_{(n+1,p)}=T^{\dbU_n(\dbZ/p)}(G)$ holds for every $n\geq2$.
By Remark \ref{G and G(p)} one can replace here $G=G_F$ by its maximal pro-$p$ Galois group $G_F(p)$.
This conjecture  known in the following cases:
\begin{enumerate}
\item[(1)]
$n=2$  (Theorem \ref{Efrat Minac Spira}(a)).
\item[(2)]
The maximal pro-$p$ quotient $G(p)$ of $G$ is a free pro-$p$ group (\cite{Efrat14a}*{Th.\ A}, \cite{MinacTan15}*{Th.\ 2.7(a)}, \cite{Efrat14b}).
\item[(3)]
$F$ is $p$-rigid and $p>2$ \cite{MinacTan15}*{Th.\ 8.1}.
\end{enumerate}

For $G=G_F$ as above, Theorem \ref{equivalence for the Zassenhaus filtration} cannot be applied to prove that $G_{(n+1,p)}=T^{\dbU_n(\dbZ/p)}$ in any other situation.
Indeed, the theorem assumes that the kernel $N$ of the projection $S\to G$ is contained in $S_{(n,p)}$, so unless (1) holds, $N\leq S_{(3,p)}$.
However we have:

\begin{prop}
Let $G=G_F$ be the absolute Galois group of a field $F$ containing a root of unity of order $p$.
Let $S$ be a free profinite group, let $N$ be a closed normal subgroup of $S$ contained in $S_{(3,p)}$, and suppose that $S/N\isom G$.
Then $G(p)$ is a free pro-$p$ group.
\end{prop}
\begin{proof}
The projection $S\to G$ induces an isomorphism $S/S_{(3,p)}\isom G/G_{(3,p)}$.
By \cite{EfratMinac17}*{Th.\ A} (which strengthens \cite{CheboluEfratMinac12}),   $H^2(G/G_{(3,p)})_{2-{\rm Massey}}\isom H^2(G)$ via inflation.
Since $S$ is also realizable as an absolute Galois group \cite{FriedJarden08}*{Cor.\ 23.1.2} and $H^2(S)=0$, we have in particular $H^2(S/S_{(3,p)})_{2-{\rm Massey}}=0$.
Therefore $H^2(G)=0$, whence $H^2(G(p))=0$ \cite{CheboluEfratMinac12}*{Lemma 6.5}.
Thus $G(p)$ is a free pro-$p$ group \cite{NeukirchSchmidtWingberg}*{Cor.\ 3.5.7}.
\end{proof}

\section{Counterexamples}
\label{section on counterexamples}
The appendix of \cite{MinacTan15}, by Min\'a\v c, T\^an and the author, constructs a profinite group $G$ such that $G_{(3,p)}\neq T^{\dbU_2(\dbZ/p)}(G)$.
The following example, which uses the more recent results of \cite{Efrat17} and \cite{Efrat20a}, extends this construction, and provides (under some assumptions on $\Omega$) closed normal subgroups $N$ of free profinite groups $S$ such that $N\leq \bar T(S)$ and  $T(G)\not\leq \pi[T(S)]$, where $\pi\colon S\to G=S/N$ is the natural projection.
This shows that the equivalent conditions of Theorem \ref{The transfer theorem general form} need not hold in general.
In the special case where $\Omega=\{\omega_2\}$ (with notation as in \S\ref{section on p Zassenhaus filtration}) this recovers the above mentioned example from \cite{MinacTan15}.

\begin{exam}
\rm
Let $S$ be the free profinite group on the basis $X=\{x_1\nek x_k\}$ of $k$ elements.
We totally order $X$ by setting $x_1<\cdots <x_k$.
A word in the alphabet $X$ is a \textsl{Lyndon word} if it is lexicographically smaller than all its proper suffixes.
Fix $n\geq2$.
For a word $w=(a_1\cdots a_n)$ of length $n$ in the alphabet $X$, consider the iterated commutator of length $n$ in $S$
\[
\tau_w= [a_1,[a_2,\cdots[a_{n-2},[a_{n-1},a_n]]\cdots]].
\]
Let $\Omega$ be a set of central extensions as in \S\ref{section on T bar T} with $m=p$ prime.
Thus $\bar T(S)/T(S)$ is an $\dbF_p$-linear space (Lemma \ref{exponent}).
We assume:
\begin{enumerate}
\item[(1)]
The iterated commutators $\tau_w$ for Lyndon words $w$ of length $n$ in $X$ belong to $\bar T(S)$, and are $\dbF_p$-linearly independent in $\bar T(S)/T(S)$;
\item[(2)]
One has $k\geq |\dbU_\omega|+n-1$ for every $\omega\in\Omega$.
\end{enumerate}

For every $n-1\leq i<j\leq k$, the word $w_{ij}=(x_1x_2\cdots x_{n-2}x_ix_j)$ is Lyndon.
Set  $\tau_{ij}=\tau_{w_{ij}}$.
In view of (1), the cosets of the $\tau_{ij}$ form a basis of a linear subspace $W$ of $\bar T(S)/T(S)$.
Let  $\varphi\colon W\to\dbZ/p$ be the linear map which maps each such coset to $1\,(\hbox{mod}\,p)$.
Let $N$ be the closed normal subgroup of $S$ generated by $\tau_{n-1,n}\tau_{ij}\inv$ for all $n-1\leq i<j\leq k$.
Thus $N\leq \bar T(S)$, and $NT(S)/T(S)\leq\Ker(\varphi)$, so $\tau_{n-1,n}\not\in NT(S)$.

On the other hand, let $\omega\in\Omega$ and let $\rho\colon S/N\to\dbU_\omega$ be a continuous homomorphism.
In light of (2), the pigeonhole principle yields $n-1\leq i<j\leq k$ such that $\rho(x_iN)=\rho(x_jN)$.
Then $\rho(\tau_{n-1,n}N)=\rho(\tau_{ij}N)=1$.
This shows that $\tau_{n-1,n}N\in T(S/N)$.
Consequently,
\[
T(S/N)\not\leq NT(S)/N,
\]
i.e., $T(G)\not\leq\pi[T(S)]$ for the projection $\pi\colon S\to G=S/N$.

Specifically, in the $p$-Zassenhaus context, where $\bar T(S)=S_{(n,p)}$ and $T(S)=S_{(n+1,p)}$, (1) holds by \cite{Efrat20b}*{Prop.\ 4.4}.
In the lower $p$-central context, where $\bar T(S)=S^{(n,p)}$ and $T(S)=S^{(n+1,p)}$, (1) holds by \cite{Efrat17}*{Th.\ 8.5(b)}.

For (2) we simply take $k$ sufficiently large.
\end{exam}

Our second example refines \cite{EfratMinac11}*{Example 13.5}:

\begin{exam}
\rm
Theorem \ref{Efrat Minac Spira}(b), and therefore also Theorem \ref{MS}(b), do not extend to primes $p>2$.
Namely, there exist absolute Galois groups $G=G_F$ as above such that $G^{(3,p)}\neq \bigcap_{s=1}^2T^{\dbU_s(\dbZ/p^{3-s})}(G)$, i.e.,
\begin{equation}
\label{counterexample}
G^{(3,p)}\neq T^{\dbZ/p^2}(G)\cap T^{\dbU_2(\dbZ/p)}(G).
\end{equation}

We first claim that  (\ref{counterexample}) holds whenever $G$ is a profinite group such that $G/G^{(3,p)}$ is non-abelian and is generated by two elements.
Indeed, since $\dbU_2(\dbZ/p)^{(3,p)}=1$ and $\dbU_2(\dbZ/p)$ is not generated by two elements, it is not a quotient of $G$.
Further, all subgroups of $\dbZ/p^2$ and $\dbU_2(\dbZ/p)$, except $\dbU_2(\dbZ/p)$, are abelian.
Therefore $G/(T^{\dbZ/p^2}(G)\cap T^{\dbU_2(\dbZ/p)}(G))$ is abelian, whence our claim.

Now let $p>2$.
\cite{EfratMinac11}*{Example 13.5} gives a field $F$ containing a root of unity of order $p$, such that $G=G_F$ satisfies
\[
G/G^{(3,p)}\isom\langle\dbZ/p^2\rangle\rtimes\langle\dbZ/p^2\rangle=\langle\tilde\tau\rangle\rtimes\langle\tilde\sig\rangle,
\]
where the generators $\tilde\tau,\tilde\sig$ satisfy $\tilde\sig\tilde\tau\tilde\sig\inv=\tilde\tau^{1+p}$.
Therefore (\ref{counterexample}) holds.
\end{exam}


\begin{bibdiv}
\begin{biblist}


\bib{ChapmanEfrat16}{article}{
author={Chapman, Michael},
author={Efrat, Ido},
title={Filtrations of the free group arising from the lower central series},
journal={J.\ Group Theory},
status={a special issue in memory of O.\ Melnikov},
volume={19},
date={2016},
pages={405\ndash433},
}

\bib{CheboluEfratMinac12}{article}{
   author={Chebolu, Sunil K.},
   author={Efrat, Ido},
   author={Min{\'a}{\v{c}}, J{\'a}n},
   title={Quotients of absolute Galois groups which determine the entire Galois cohomology},
   journal={Math. Ann.},
   volume={352},
   date={2012},
   pages={205--221},
}


\bib{Dwyer75}{article}{
   author={Dwyer, William G.},
   title={Homology, Massey products and maps between groups},
   journal={J. Pure Appl. Algebra},
   volume={6},
   date={1975},
   pages={177--190},
}


\bib{Efrat14a}{article}{
author={Efrat, Ido},
title={The Zassenhaus filtration, Massey products, and representations of profinite groups},
journal={Adv.\ Math.},
volume={263},
date={2014},
pages={389\ndash411},
}

\bib{Efrat14b}{article}{
author={Efrat, Ido},
title={Filtrations of free groups as intersections},
journal={Archiv d.\ Math},
volume={103},
date={2014},
pages={411\ndash420},
}

\bib{Efrat17}{article}{
author={Efrat, Ido},
title={The cohomology of canonical quotients of free groups and Lyndon words},
journal={Documenta Math.},
volume={22},
date={2017},
pages={973\ndash997},
}

\bib{Efrat20a}{article}{
author={Efrat, Ido},
title={The lower p-central series of a free profinite group and the shuffle algebra},
journal={J.\ Pure Appl.\ Algebra},
volume={224},
date={2020},
pages={106260}
}

\bib{Efrat20b}{article}{
author={Efrat, Ido},
title={The p-Zassenhaus filtration of a free profinite group and shuffle relations},
date={2020},
status={to appear},
eprint={arXiv:2003.08903}
}


\bib{EfratMinac11}{article}{
label={EfMi11},
   author={Efrat, Ido},
   author={Min\'a\v c, J\'an},
   title={On the descending central sequence of absolute Galois groups},
   journal={Amer. J. Math.},
   volume={133},
   date={2011},
   pages={1503\ndash1532},
 }


\bib{EfratMinac17}{article}{
label={EfMi17},
   author={Efrat, Ido},
   author={Min\'a\v c, J\'an},
   title={Galois groups and cohomological functors},
   journal={Trans.\ Amer.\ Math.\ Soc.},
   volume={369},
   date={2017},
   pages={2697\ndash2720},
 }

\bib{FriedJarden08}{book}{
   author={Fried, Michael D.},
   author={Jarden, Moshe},
   title={Field Arithmetic},
   edition={Third edition},
   publisher={Springer, Berlin},
   date={2008},
   pages={xxiv+792},
}



\bib{GilleSzamuely06}{book}{
 author={Gille, Philippe},
   author={Szamuely, Tam{\'a}s},
   title={Central Simple Algebras and Galois Cohomology},
   publisher={Cambridge University Press},
   place={Cambridge},
   date={2006},
   pages={xii+343},
}

\bib{Hoechsmann68}{article}{
author={Hoechsmann, Klaus},
title={Zum Einbettungsproblem},
journal={J.\ reine angew.\ Math.},
volume={229},
date={1968},
pages={81\ndash106},
}


\bib{MerkurjevSuslin82}{article}{
    author={Merkurjev, A. S.},
    author={Suslin, A. A.},
     title={$K$-cohomology of Severi-Brauer varieties and the norm residue homomorphism},
  language={Russian},
   journal={Izv. Akad. Nauk SSSR Ser. Mat.},
    volume={46},
      date={1982},
     pages={1011\ndash 1046},
    translation={
        journal={Math. USSR Izv.},
         volume={21},
           date={1983},
   pages={307\ndash 340},
} }



\bib{MinacSpira96}{article}{
  author={Min{\'a}{\v {c}}, J{\'a}n},
  author={Spira, Michel},
  title={Witt rings and Galois groups},
  journal={Ann. Math.},
  volume={144},
  date={1996},
  pages={35\ndash60},
  label={MSp96},
}


\bib{MinacTan15}{article}{
   author={Min\'a\v c, J\'an},
   author={T\^an, Nguyn Duy},
   title={The kernel unipotent conjecture and the vanishing of Massey    products for odd rigid fields},
   journal={Adv. Math.},
   volume={273},
   date={2015},
   pages={242--270},
   status={(with an appendix by I.\ Efrat, J.\ Min\'a\v c, and N.D. T\^an)},
}


\bib{NeukirchSchmidtWingberg}{book}{
author={Neukirch, J.},
author={Schmidt, Alexander},
author={Wingberg, Kay},
title={Cohomology of Number Fields},
edition={Second edition},
publisher={Springer},
date={2008},
}

\end{biblist}
\end{bibdiv}
\end{document}